\PassOptionsToPackage{colorlinks=true, allcolors=blue}{hyperref}  %
\documentclass[a4paper,fleqn]{cas-dc}

\usepackage{float}
\usepackage{amsthm}
\usepackage{mdframed}
\usepackage{amsthm}
\usepackage{siunitx}
\usepackage{graphicx}
\usepackage{dsfont}
\usepackage[ruled]{algorithm2e}
\SetKwInput{Input}{Input}
\usepackage{braket}
\usepackage{subcaption}  %
\usepackage{mathrsfs}
\usepackage{accents} %
\usepackage[numbers,sort&compress]{natbib}

\newtheorem{theorem}{Theorem}

\def\tsc#1{\csdef{#1}{\textsc{\lowercase{#1}}\xspace}}
\tsc{WGM}
\tsc{QE}

\begin{document}
\let\WriteBookmarks\relax
\def\floatpagepagefraction{1}
\def\textpagefraction{.001}

\shorttitle{}

\shortauthors{}
\ExplSyntaxOn
\cs_set:Npn \__first_footerline: { }
\ExplSyntaxOff

\title [mode = title]{Quantum-Inspired Tensor Train Spectral Solver for High-Resolution PDEs with Non-Periodic Boundary Conditions}

\author[1, 2]{Sascha H. Hauck}[orcid=0009-0002-4800-0598]%
\cormark[1]
\ead{sascha.hannes.hauck@itwm.fraunhofer.de}
\author[3]{Nis-Luca van Hülst}%

\author[1]{Matthias Kabel}%

\author[3,4]{Dieter Jaksch}%

\author[2]{Nicolas R. Gauger}%

\affiliation[1]{organization={Fraunhofer ITWM},
            addressline={Department of Flow and Material Simulation},
            city={Kaiserslautern},
            postcode={67663},
            state={Rhineland Palatinate},
            country={Germany}}

\affiliation[2]{organization={University of Kaiserslautern-Landau (RPTU)},
            addressline={Chair for Scientific Computing},
            city={Kaiserslautern},
            postcode={67663},
            state={Rhineland Palatinate},
            country={Germany}}

\affiliation[3]{organization={University of Hamburg},
            addressline={Institute for Quantum Physics},
            city={Hamburg},
            postcode={20146},
            state={Hamburg},
            country={Germany}}

\affiliation[4]{organization={University of Oxford}, addressline={Clarendon Laboratory, Parks Road}, city={Oxford}, postcode={OX1 3PU}, country={United Kingdom}}

\cortext[1]{Corresponding author}

\begin{abstract}
Tensor train (TT) methods enable high-resolution computational mechanics simulations by exploiting low-rank structure in discretized fields and operators.
While the superfast Fourier transform (SFFT) has enabled efficient TT-based spectral solvers for periodic problems, corresponding representations for non-periodic boundary conditions remain underexplored.
We introduce the superfast cosine and sine transforms (SFCT and SFST), extending the TT spectral framework to homogeneous Neumann and Dirichlet boundary conditions, respectively.
Building on a quantum-circuit construction of the discrete cosine and sine transforms, we derive both transforms from a shared TT operator representation and prove that their TT ranks are bounded by \(6r_S\), independently of the discretization order, where \(r_S\) is the rank of the SFFT at the prescribed accuracy.
The measured ranks plateau near 32 for grid sizes up to \(2^{26}\) points.
For Rayleigh--B\'enard convection, spectral pressure projection
achieves up to a 7-fold speedup in the Poisson solve relative to a
preconditioned physical-space TT solver while maintaining comparable accuracy.
Three-dimensional computational homogenization achieves a 2--6-fold
speedup and a 4--7-fold reduction in memory consumption compared with
the existing TT-based approach.
The results establish an efficient TT framework for spectral PDE solvers with non-periodic boundary conditions.
\end{abstract}


\ExplSyntaxOn\keys_set:nn { stm / mktitle } { nologo = true }\ExplSyntaxOff
\maketitle

\section{Introduction}

Partial differential equations (PDEs) model a wide range of physical
phenomena, including fluid flow, heat transport, and the response of
heterogeneous materials. Their numerical solution
on fine grids can require substantial memory and computing time.
Tensor-Train (TT) methods approximate the discretized solution by a chain
of small tensors~\cite{Oseledets2011,Schollwock2011}. The dimensions
connecting these tensors, called TT ranks, determine the storage cost
and the attainable approximation accuracy. Calculations can be carried
out on the compressed representation, with potential savings when the
ranks remain moderate. TT methods have proven beneficial in a wide range of high-resolution applications, including fluid dynamics~\cite{Gourianov2022,Kiffner2023}, kinetic plasma dynamics~\cite{Ye2022}, simulations of the Gross--Pitaevskii equation~\cite{ConnorGP,GomezGP}, structural mechanics problems~\cite{Ali2025,Risthaus2023} and computational homogenization~\cite{Vondrejc2020,Hauck2026,Hauck2027, Risthaus2024}.

For grid dimensions of the form $2^n$ with $n\in\mathds{N}$, TT representations provide a link between tensor-based numerical methods and quantum computing (see~\cite{CecReview} for a review).
This connection can be exploited in both directions: tensor networks have been used to construct quantum circuits for differential equations~\cite{SieglPRR} and hybrid algorithms for nonlinear flow problems~\cite{SieglHybrid}, while quantum-circuit constructions can, conversely, be used to derive low-rank TT representations.
For the quantum Fourier transform (QFT), such a low-rank representation has been established in the form of the Superfast Fourier Transform (SFFT)~\cite{Dolgov2012,Chen2023}. The SFFT applies the transform directly to fields in compressed TT form, avoiding the reconstruction of the full grid vector. Its computational efficiency depends on the intermediate TT ranks. For periodic problems, this approach enables direct Fourier-space operations within a TT solver~\cite{Hauck2026}.

While the SFFT provides an efficient approach for periodic problems, practical applications also involve non-periodic boundary conditions, such as Dirichlet and Neumann conditions, or combinations thereof. Spectral methods offer a natural way to accommodate these conditions through basis functions adapted to the underlying boundary-value problem. For suitable operators on uniform Cartesian grids, discrete cosine and sine transforms diagonalize the corresponding finite-difference operators for homogeneous Neumann and Dirichlet conditions, respectively. This allows the discrete equations to be solved in spectral space, although the appropriate transform depends on the grid arrangement and boundary discretization~\cite{Strang1999}. Such operators arise, for example, in pressure-correction equations for wall-bounded incompressible flow~\cite{Guermond2006}
and in computational homogenization with prescribed displacement or traction boundary conditions~\cite{Grimm2021,Javili2017}.

TT methods for non-periodic boundary-value problems include finite-element
solvers~\cite{Marcati2022,Kornev2024} and explicit representations of the
one-dimensional discrete inverse Laplacian~\cite{Kazeev2012}.
For thermal homogenization, Risthaus and Schneider~\cite{Risthaus2024}
compared a solver based on conventional sine and cosine transforms with one
based on a QTT approximation of the Green's operator.
However, a direct TT representation of the spectral transforms
themselves has not yet been established. To address this gap, we
adapt the quantum-circuit construction of~\cite{Klappenecker01} and derive the
type-II superfast cosine and sine transforms (SFCT and
SFST) from a shared QFT-based TT operator. We prove that their ranks
are bounded by \(6r_S\), where \(r_S\) is the SFFT rank,
independently of the grid size.
We evaluate their numerical accuracy, ranks,
storage requirements, and runtime for sizes up to \(2^{26}\) grid points.
Finally, we demonstrate the applicability of the transforms to two relevant problems in computational mechanics.
For Rayleigh--B\'enard convection, we compare a mixed
Fourier--cosine pressure solver with a preconditioned physical-space TT
solve. In three-dimensional computational homogenization, we apply the
spectral transforms to problems with mixed, periodicity-compatible uniform
boundary conditions.

The remainder of this paper is organized as follows.
Section~\ref{sec:tensor_trains} summarizes the basics of the TT format.
In Section~\ref{sec:spectral_methods}, we develop the SFCT and SFST and use their quantum-circuit construction to establish the expressivity theorem. We then evaluate their accuracy and computational cost numerically.
Sections~\ref{sec:quantum_inspired_cfd} and~\ref{sec:quantum_inspired_homogenization} present two complementary applications in computational mechanics: Rayleigh--B\'enard convection and computational homogenization.
Finally, Section~\ref{sec:conclusion} summarizes the main developments and findings of this work.

\section{Tensor-Train Format}
\label{sec:tensor_trains}
In this work, we employ the tensor train (TT) format to represent a tensor $T$ of order $n$ as a chain of interconnected third-order tensors, referred to as TT cores.
The decomposition of $T$ is given by~\cite{Oseledets2011}
\begin{align}
\label{eq:tt_format}
T(i_1, ..., i_n) &= \sum_{\boldsymbol{\alpha}}
A^{\alpha_1}_{\alpha_0}(i_1) \cdots A^{\alpha_n}_{\alpha_{n-1}}(i_n) \\
&= A_1(i_1) \cdots A_n(i_n),
\end{align}
where $A_k$ denotes the $k$-th TT core with physical dimension $d$ and shape $(r_{k-1}\, \times \, d \, \times \, r_k)$.
The tuple $\boldsymbol{\alpha}=(\alpha_0, \ldots, \alpha_k, \ldots, \alpha_n)$ collects the virtual indices.
The corresponding TT ranks $r_k$ determine the expressive capacity of the representation and, consequently, its computational cost.
The boundary ranks are conventionally fixed to the trivial dimension $r_0=r_n=1$ to obtain a closed representation.

For an $n$-th order tensor with uniform physical dimension $d$, direct storage requires $O(d^n)$ elements, whereas the TT representation requires only $O(ndr^2)$ parameters, provided the TT ranks are bounded by $r$.
Thus, for sufficiently low ranks,
the representation reduces the storage requirement from exponential to linear in the tensor order $n$.
This low-rank structure also enables efficient manipulation of high-dimensional tensors without explicitly forming the full tensor \cite{Hackbusch2009, Hackbusch2015}, potentially leading to significant runtime improvements in the process.

To apply the TT format to spatially discretized fields,
such as those considered throughout this work,
we first reshape their sampled values into high-order
tensors by expressing the grid indices in binary form.
For a uniform grid with $2^n$ points on $[0,1)$, each coordinate
can be expressed through
\begin{align}
\label{eq:binary_index}
    x(i_1,\ldots,i_n)
    = \sum_{m=1}^{n} i_m 2^{-m},
    \qquad i_m\in\{0,1\}.
\end{align}
Each grid point is therefore identified by a binary string
$(i_1,\ldots,i_n)\in\{0,1\}^{n}$.
Generalizing this idea to a three-dimensional field sampled on a tensor-product grid results in an order-$3n$ tensor.
Thus, defining the tuples $\boldsymbol{i}=(i_1,\ldots,i_n)$ and analogously
$\boldsymbol{j}$ and $\boldsymbol{k}$, we define the high-order tensor
\begin{align*}
    F(\boldsymbol{i},\boldsymbol{j},\boldsymbol{k})
    := f\bigl(
        x(\boldsymbol{i}),
        y(\boldsymbol{j}),
        z(\boldsymbol{k})
    \bigr).
\end{align*}
Since this procedure is solely a bookkeeping operation, it introduces no additional approximation to the sampled field.
The resulting high-order tensor can be transferred into TT format by employing the TT-SVD algorithm~\cite{Oseledets2011}, yielding the
quantized tensor-train (QTT) representation, with physical dimension
$d=2$. This representation is mathematically analogous to the tensor-product structure of a quantum many-body system, where each binary mode corresponds to a single qubit with a two-dimensional Hilbert space \cite{Fannes1992}.
This analogy will be exploited in subsequent sections.

Scalar multiplication, addition, element-wise multiplication, and the application of linear operators can be performed directly in TT format. Scalar multiplication leaves the ranks unchanged, whereas addition, element-wise multiplication, and operator application can increase the TT ranks, potentially leading to an increasingly costly representation as successive operations are performed.
However, TT rounding reduces these ranks
while controlling the relative Frobenius error $\varepsilon_\mathrm{TT}$~\cite{Oseledets2011}.
Alternatively, the TT ranks can be truncated to a prescribed maximum value, providing direct control over the representation complexity at the expense of potentially larger approximation errors.

For the application of a general linear map, the operator itself must also be represented in a suitable tensor format. A tensor-train operator (TTO) represents the linear map in TT form and can be applied to a TT tensor through local contractions between the corresponding operator and tensor cores. In analogy to the TT representation, each TTO core contains two physical indices, corresponding to the input and output spaces. Therefore, a TTO is written as
\begin{align}
    \label{eq:tto_format}
    \mathcal{T}(i_1,j_1,\ldots,i_n,j_n)
    &=
    \mathcal{G}_1(i_1,j_1)
    \cdots
    \mathcal{G}_n(i_n,j_n),
\end{align}
where $\mathcal{G}_k$ denotes the $k$-th TTO core.
For the square operators considered here, each TTO core has dimensions
$(\rho_{k-1}\times d\times d\times\rho_k)$, where $\rho_k$ denotes
the operator rank and $\rho_0=\rho_n=1$.
The indices $i_k$ and $j_k$ correspond to the output and input,
respectively. Applying a TTO to a TT contracts each input index
$j_k$ with the corresponding physical index of the TT.
Before compression, the resulting representation has bond dimensions
$\rho_k r_k$.

Another central operation is the Hadamard product $\odot$, which refers to the
pointwise multiplication of discretized fields.
Such operations arise, for example, in the evaluation of nonlinear terms and the application of diagonal operators in spectral space, and can constitute a major computational bottleneck in TT-based solvers.
For two fields $x$ and $y$ in TT format, the Hadamard product can be written as
$x\odot y=\operatorname{diag}(y)x$.
If $A_k$ is the $k$-th TT core of $y$, the corresponding core of $\operatorname{diag}(y)$ is defined as
\begin{align}
 \label{eq:tt_to_tto}
 \mathcal{G}_k(i_k,j_k)
 =
 A_k(i_k)\,\delta_{i_kj_k}.
\end{align}
Applying the diagonal TTO to $x$ is a standard TTO--TT contraction.
If both input TTs have rank $r$, the product ranks can reach $r^2$.
Reducing these ranks using standard TT rounding costs $O(nr^6)$.

Variational compression (VAR) instead fits an output TT $z$ without forming the rank-$r^2$ intermediate.
This fit can be computed via alternating least squares (ALS)~\cite{Holtz2012}, which minimizes $\|z-\operatorname{diag}(y)x\|_F^2$ over one core of $z$ at a time.
For comparable input and output ranks, each sweep costs $O(nr^4)$~\cite{Michailidis2024TTM}. Rank-adaptive alternatives include the two-site modified alternating least squares (MALS) method~\cite{Holtz2012}, which updates two neighboring cores together, and the alternating minimal energy method (AMEn)~\cite{Dolgov2013}, which enriches the local basis using residual information.

Algorithms for the Hadamard product that can achieve cubic scaling in the TT rank include tensor train
multiplication (TTM)~\cite{Michailidis2024TTM}, alternating cross
interpolation (ACI)~\cite{Ritter2026ACI}, and recursive sketched
interpolation (RSI)~\cite{Meng2026-arxiv}.
In Section~\ref{sec:quantum_inspired_cfd},
we test RSI as an alternative to VAR for the Hadamard product within the spectral Poisson solver.

\section{Boundary-Adapted Spectral Transforms in TT Format}
\label{sec:spectral_methods}

In this section, we develop TT representations of boundary-adapted spectral transforms, connecting classical spectral methods with quantum-circuit constructions and their classical tensor-network representations. We focus on the type-II discrete cosine and sine transforms (DCT-II and DST-II), which are relevant to the boundary conditions considered in this work.
The construction can be extended analogously to other types of cosine and sine transforms.
Starting from the classical construction of these transforms, we derive a unitary quantum formulation and exploit its circuit structure to establish an upper bound on the TT ranks of the resulting transform operators. This analysis provides the foundation for the Superfast Cosine and Sine Transforms introduced in the subsequent sections.

\subsection{Classical Spectral Methods}

To illustrate the construction of spectral methods for boundary value problems, we briefly review Fourier-type transforms for solving elliptic partial differential equations on bounded domains. As a simple example, we consider the one-dimensional Poisson problem to establish the fundamental concepts underlying the subsequent formulations.
Let $\Omega=(0,L)\subset \mathbb{R}$ and consider the problem
\begin{equation}
\label{eq:poisson}
- \frac{\mathrm{d}^2 u(x)}{\mathrm{d}x^2} = f(x), \qquad \forall x \in \Omega,
\end{equation}
subject to homogeneous Neumann boundary conditions
\begin{align*}
    \frac{\mathrm{d} u(x)}{\mathrm{d}x}\big|_{x=0} =\frac{\mathrm{d} u(x)} {\mathrm{d}x}\big|_{x=L} =0.
\end{align*}
In general, the continuous problem is reformulated in a discretized setting to obtain a computationally tractable formulation.
Therefore, the domain $\Omega$ will be represented using a cell-centered grid
\[
x_m=\left(m+\frac12\right)h,
\qquad
h=\frac{L}{N},
\]
with $m=0,\ldots,N-1$.
For consistency with the subsequent sections, we assume that the number of grid points is a power of two, $N=2^n$, with $n\in\mathbb{N}$.
The continuous functions in Eq.~\eqref{eq:poisson} are sampled at the grid points, yielding the discrete values $u_m=u(x_m)$ and $f_m=f(x_m)$ for the field and source terms, respectively.
We then construct the even extension of $u$ as
\begin{align*}
\widetilde{u}=[u_0,u_1,\cdots,u_{N-1},u_{N-1},\cdots,u_1,u_0]^T,
\end{align*}
which contains $2N$ entries and is indexed similarly by $m=0,\ldots,2N-1$.
An analogous even extension is constructed for $f$.

The discrete Fourier transform (DFT) can then be applied to these extended vectors, yielding
\begin{align}
\label{eq:u_fourier}
\widetilde{u}_m &=
\frac{1}{\sqrt{2N}} \sum_{k=0}^{2N-1} \hat{u}_{k} e^{i \pi \frac{m k}{N}}, \\
\label{eq:f_fourier}
\widetilde{f}_m &=
\frac{1}{\sqrt{2N}} \sum_{k=0}^{2N-1} \hat{f}_{k} e^{i \pi\frac{ m k}{N}},
\end{align}
with the Fourier coefficients $\hat{u}_{k}$ and $\hat{f}_{k}$, respectively.
Additionally, the component wise Poisson equation \eqref{eq:poisson} can be written through
\begin{align}
\label{eq:poisson_disc}
    (\mathcal{L}u)_m \equiv \frac{u_{m+1}-2 u_m + u_{m-1}}{h^2} = -f_m.
\end{align}
To close the above finite difference operation at the boundary points, we introduce two ghost points through the discrete Neumann conditions $u_{-1}=u_0$ and $u_N=u_{N-1}$.

Substituting the transformations \eqref{eq:u_fourier} and \eqref{eq:f_fourier} into the component-wise discretized Poisson equation \eqref{eq:poisson_disc} yields an alternative representation of the original problem in Fourier space
\begin{align}
    \sum_{k=0}^{2N-1} e^\frac{i \pi k m}{N}\left[ \frac{\hat{u}_k}{h^2} \left(e^\frac{i \pi k}{N} + e^\frac{-i \pi k}{N} -2\right) + \hat{f}_k\right] = 0.
\end{align}
The differential equation we started with is thus transformed into an algebraic equation in Fourier space.
Owing to the orthogonality of the Fourier modes, the above expression further reduces to the component-wise equation
\begin{align*}
\xi_k^2 \hat{u}_{k} = \hat{f}_{k},
\end{align*}
with $\xi_k = 2 N \sin\left[\pi k / (2N)\right]/ L$ for $k\in\{1, \cdots, 2N-1\}$.
Since $\xi_0=0$, the zero mode $k=0$ imposes the solvability condition that discrete source has zero mean
\begin{align}
    \label{eq:poisson_compatibility}
    \hat{f}_0 =
    \frac{1}{\sqrt{2N}}\sum_{m=0}^{2N-1} \widetilde{f}_m = 0.
\end{align}
Moreover, the zero mode implies that the solution $u$ is not unique, since its mean component remains undetermined.
In particular, if $u$ is a solution, then so is $u\rightarrow u + C$ for any constant $C$.

The even extension allows the Poisson problem to be solved using the DFT, but doubles the number of grid points from \(N\) to \(2N\) without adding new information. The discrete cosine transform avoids this redundancy by incorporating the even symmetry directly into the transform, allowing the original \(N\)-point problem to be solved without explicitly constructing the extended signal.

To see how this is achieved, let us consider the Fourier transform of the even extension $\widetilde{u}$
\begin{align*}
   \hat{u}_k
   &=
   \frac{1}{\sqrt{2N}} \sum_{m=0}^{2N-1} \widetilde{u}_m e^{-i \pi \frac{k m}{N}}
   \\
    &=
   \frac{1}{\sqrt{2N}} \sum_{m=0}^{N-1} \left[\widetilde{u}_m e^{-i \pi \frac{ k m}{N}} + \widetilde{u}_{2N-1-m} e^{-i \pi\frac{k (2N-1-m)}{N}}\right]
   \\
    &=
   \frac{e^{i \frac{\pi k}{2N}} }{\sqrt{2N}} \sum_{m=0}^{N-1} \left[\widetilde{u}_m e^{-i\pi \frac{ k (2m+1)}{2N}} + \widetilde{u}_{2N-1-m} e^{i \pi\frac{k (2m+1)}{2N}}\right].
\end{align*}
Due to the even extension, the symmetry condition $\widetilde{u}_m= \widetilde{u}_{2N-1-m}\equiv u_m$ is fulfilled.
Therefore we are able to simplify the above relation to
\begin{align*}
   \hat{u}_k
   &=
   \sqrt{\frac{2}{N}} e^{i \frac{\pi k}{2N}}\sum_{m=0}^{N-1} u_m \cos\left[\frac{\pi k}{N} \left(m+\frac{1}{2}\right)\right] \\
   &\equiv e^{i \frac{\pi k}{2N}} b_k^{-1}C^\mathrm{II}_N[u]_k,
\end{align*}
with the orthonormalization factor $b_k=1$ for $k\neq 0$ and $b_0=1/\sqrt{2}$.

The above expression connects the original $N$-point vector $u$ with its components in Fourier space via the orthonormal DCT-II transform $C_N^{\mathrm{II}}$. Therefore, we can solve the one-dimensional Poisson problem with Neumann conditions without explicit even extension.

Similar derivations can be performed for different discretization schemes, leading to the different DCT types, as well as for Dirichlet boundary conditions, leading to the discrete sine transform (DST) \cite{Wickerhauser93, Strang1999}.
More complex inhomogeneous Dirichlet boundary conditions can be incorporated via a lifting function $g(x)$ satisfying the prescribed boundary values. Writing $u=v+g$ yields a homogeneous problem for $v$ with a modified right-hand side. Neumann and mixed boundary conditions can be treated analogously by choosing appropriate symmetric or antisymmetric extensions, leading to cosine- or mixed-sine spectral bases \cite{dePrenter2023}.

\subsection{Quantum Spectral Methods}
\label{subsec:quantum_spectral_methods}
Quantum Spectral Methods build on variants of the Quantum Fourier Transform (QFT) adapted to different boundary conditions.
While periodic boundary conditions are treated using algorithms based directly on the QFT, Dirichlet and Neumann boundary conditions are naturally associated with the Quantum Sine Transform (QST) and Quantum Cosine Transform (QCT), respectively.

\begin{figure*}
    \centering
    \includegraphics[width=\linewidth]{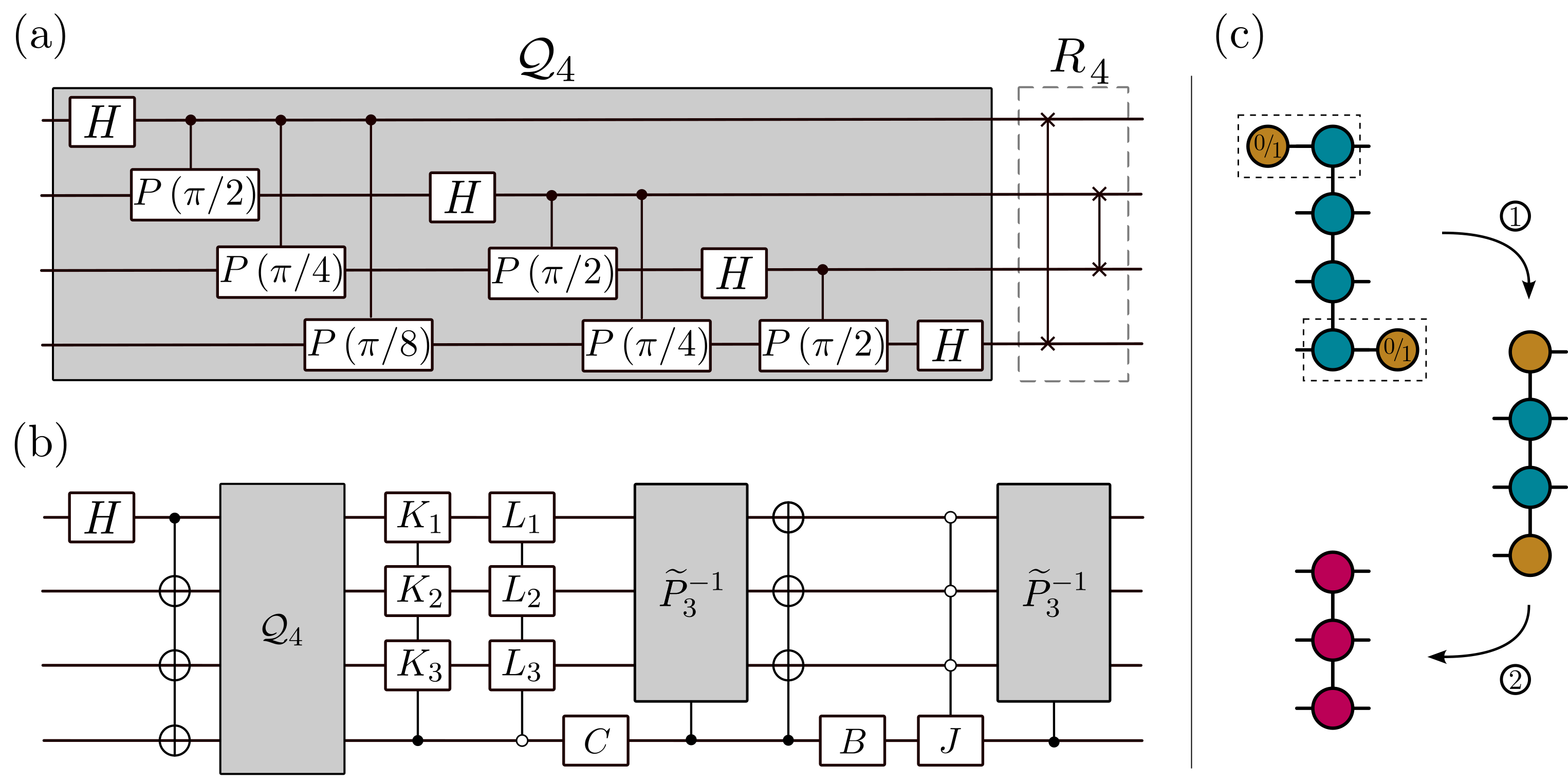}
    \caption{
    Quantum circuit for $n=4$ qubits:
    (a) the QFT and
    (b) the QCST.
    (c) Extraction of either the SFCT or SFST from their shared TTO representation.
    Selecting either SFCT or SFST is achieved by fixing the boundary cores to trivial $0-$ or $1-$vectors. This procedure generally produces a skewed TTO representation, which is subsequently restored to standard canonical TTO form via a reshaping step.
    The circuits in (a) and (b) are based on the constructions presented in \cite{Nielsen2010} and \cite{Klappenecker01}, with modifications described in the text.
    }
    \label{fig:qft_qc}
\end{figure*}

The quantum circuit implementation of the QFT is shown in Figure~\ref{fig:qft_qc} (a), where the $n$-qubit QFT is represented through a product
\begin{align}
\label{eq:reduced_qft_form}
    \mathcal{QFT}_n=\mathcal{R}_n \mathcal{Q}_n
\end{align}
of the pure transformation $\mathcal{Q}_n$ and a reordering operation $\mathcal{R}_n$~\cite{Nielsen2010}.
While the pure transformation depends solely on elementary gates, such as the Hadamard $H$ and the phase gate $P(\theta)$, the reordering operation depends on multiple SWAP gates.
These rely heavily on CNOT gates, which remain a comparatively challenging component even in the quantum setting \cite{Zheng2024}.
Moreover, the reordering is not necessarily needed for the actual transform but rather is introduced in the literature to match the ordering conventions used in classical DFT methods.
In contrast to its classical counterpart, which exhibits an exponential time-complexity in the number of qubits, requiring $\mathcal{O}(n 2^n)$ operations, the QFT achieves a quadratic time complexity of $\mathcal{O}(n^2)$ \cite{Coppersmith1994, Nielsen2010}.

Classically, the DCT and DST can be derived from the DFT using even and odd extensions of the input data, respectively. However, simply discarding redundant Fourier coefficients is not a unitary operation and therefore cannot be directly implemented on a quantum computer.
A suitable quantum circuit must implement the desired transform while preserving unitarity.

This obstacle can be avoided by embedding the sine and cosine transforms into a single unitary framework rather than treating them as independent reductions of the Fourier transform. Specifically, one can decompose the Fourier transform $F_{2N}$ into two symmetry-adapted components acting on complementary subspaces
\begin{align}
\label{eq:SpectralRelation}
    U^\dagger_N F_{2N} V_N = C^\mathrm{II}_N \oplus (-i) S^\mathrm{II}_N,
\end{align}
where $C^\mathrm{II}_N$ and $S^\mathrm{II}_N$ denote the discrete cosine and sine transforms of type II, respectively \cite{Wickerhauser93}.
Similar decompositions can be derived for other types of discrete cosine and sine transforms.
This identity defines a unitary mapping on the full Hilbert space of dimension $2N=2^{n+1}$.
It therefore admits a natural quantum generalization by replacing the classical FFT with the QFT and constructing quantum circuits for the auxiliary basis transformations $U_N$ and $V_N$ \cite{Klappenecker01}.

The full quantum circuit implementing the combined Quantum Cosine-Sine Transform (QCST), analog to Eq.~\eqref{eq:SpectralRelation},
is shown in Fig.~\ref{fig:qft_qc}(b).
We emphasize that the circuit is formulated exclusively in terms of the reduced transform $\mathcal{Q}_n$,
cf. Eq.~\eqref{eq:reduced_qft_form}.
While this restriction is not strictly necessary for the quantum implementation,
it is crucial for obtaining an efficient classical tensor train representation.
The full circuit representation employs a set of non-standard gates defined as
\begin{align*}
    B = \frac{1}{\sqrt{2}}
    \begin{pmatrix}
        1 & i\\
        1 & -i
    \end{pmatrix}
    \qquad
    J = \frac{1}{\sqrt{2}}
    \begin{pmatrix}
        1 & -i\\
        -i & 1
    \end{pmatrix}
\end{align*}
as well as
$K_j = \text{diag}(\overline{\omega}^{2^{j-1}}, 1)$,
$L_j = \text{diag}(1, \omega^{2^{j-1}})$
and $C = \text{diag}(1, \overline{\omega})$
where $\omega=\exp\left[2\pi i/(4N)\right]$.
In addition, the construction makes use of the bit-reversed inverse permutation gate
\begin{align*}
\widetilde{P}_n^{-1} &= R_\text{inv} P_n^{-1} R_\text{inv}\\
    P_n^{-1}&:\ket{x} \rightarrow \ket{(x-1)\text{mod}\, 2^n},
\end{align*}
with $x\in\{0,\ldots,2^n-1\}$.
The operator $R_\text{inv}$ constitutes the bit reversal operation $R_\text{inv}\ket{b_1 b_2 \cdots b_n} = \ket{b_n \cdots b_2 b_1}$.
The original permutation gate $P_n$ is a well-known construction consisting solely of multi-controlled NOT gates \cite{Pueschel99}.

\subsection{Quantum-Inspired Spectral Methods}

Having established a full quantum circuit representation of the QCST, we now address its classical simulability. In particular, we ask how the complexity of this transformation scales in a classical TT representation.
It is well known that the reduced QFT $\mathcal{Q}_n$, which forms the backbone of the QCST construction, admits an efficient classical representation in TT format known as the Superfast Fourier Transform (SFFT) \cite{Dolgov2012, Chen2023}.
This suggests that similar properties are inherited by the QCST structure as well.
This leads to one of the central results of this work:
\begin{theorem}[Expressivity Theorem]
\label{theo:expressivity_theorem}
The quantum circuit implementing the QCST admits a TT representation with maximal TT rank $\hat{r}$ bounded by
\begin{align*}
    \hat{r} \leq 12\,r_S,
\end{align*}
where $r_S$ denotes the TT rank of the SFFT associated with a relative accuracy $\varepsilon_\mathrm{TT}$. If the circuit is restricted to implement either the quantum version of the DCT or the DST individually, the corresponding TT rank satisfies the tighter bound
\begin{align*}
    \hat{r} \leq 6\,r_S.
\end{align*}
\end{theorem}
The proof takes the operator Schmidt decomposition of the circuit in
Fig.~\ref{fig:qft_qc}(b) across every qubit bipartition, left and right of the
reduced QFT $\mathcal{Q}_n$. The initial Hadamard--CNOT block has Schmidt rank
two. To the right of $\mathcal{Q}_n$, the controlled $K/L$ block reduces to a
product operator of rank one, while the two controlled permutations have rank
three and the CNOT-$B$-$J$ block has rank four. Because these components share
operator factors, the assembled right-hand block does not reach the naive
product of these ranks but has Schmidt rank exactly six, which drops to three
once the output index of the last qubit is fixed. Combining the leading factor of two with the rank $r_S$ of $\mathcal{Q}_n$, which is determined by the prescribed operator accuracy $\varepsilon_\mathrm{TT}$, yields the bounds
$\hat{r} \leq 2 \cdot 6\,r_S$ and $\hat{r} \leq 2 \cdot 3\,r_S$, respectively.
A detailed proof of Theorem~\ref{theo:expressivity_theorem} can be found in
Appendix~\ref{app:proof_expressivity_theorem}.
Following the naming convention used for the SFFT, we refer to the TTO representations of the QCT and QST as the Superfast Cosine Transform (SFCT) and Superfast Sine Transform (SFST), respectively.

The tighter bound applies directly to the SFCT and SFST: both are obtained from the QCST by fixing the boundary physical indices, see Fig.~\ref{fig:qft_qc}(c). Since fixing an index cannot increase a bond dimension, the ranks of the transforms are bounded by those of the QCST. Reshaping the resulting tensors into standard TTO form (step~2) may introduce an additional rank increase, but the observed ranks remain well below $6\,r_S$, as shown in Sec.~\ref{subsec:numerical_results}.
The ranks required to achieve a prescribed accuracy $\varepsilon_\mathrm{TT}$ for the SFFT saturate as the discretization is refined. The rank bounds derived above therefore imply the same saturation behavior for the QCST and, consequently, for the SFCT and SFST.
The question remains under which conditions the TTO version is preferred over
the classical transforms.

Before proceeding with the analysis, we convert the quantum circuit into a TTO representation using the so-called Zip-Up algorithm, which constructs the operator sequentially, one slice at a time \cite{Stoudenmire2010}. The resulting slices are successively contracted and compressed to obtain a single TTO representing the full transformation.

Importantly, within the TN framework the unitarity constraint is no longer required. This allows us to construct separate TTO representations for the cosine and sine transforms from the full QCST operator without discarding information explicitly. To this end, we exploit the block-diagonal structure of the decomposition in Eq.~\eqref{eq:SpectralRelation}.
The DCT and DST components can be extracted by fixing the boundary physical indices of the TTO. Concretely, we project the input index of the first core and the output index of the last core onto basis states $\ket{0}$ or $\ket{1}$, see Fig.~\ref{fig:qft_qc}(c).
This selection reflects the structure of the underlying operator: the input index of the first core determines whether the upper or lower block of the matrix representation in Eq.~\eqref{eq:SpectralRelation} is selected, while the output index of the last core distinguishes between the left and right blocks. The latter arises from the adopted momentum-space convention, which reverses the ordering of the cores. Consequently, the first core encodes both the most significant input bit and the least significant output bit.

Performing the resulting contractions yields a TTO in a slightly distorted form, where the boundary cores carry either a single input or a single output index,
see Fig.~\ref{fig:qft_qc}(c), step~1.
To restore the standard TTO structure, we apply a sequence of contractions and reshaping operations.
First, the uppermost core is contracted with its neighboring core and reshaped to expose one input and two output indices. The resulting tensor is contracted with the next core and reshaped again, effectively transferring the output index to the next core.
This is followed by a SVD, which splits the compound tensor back into two cores: the first regains the standard TTO structure, while the second temporarily carries two output indices.
Repeating this procedure moves the output indices one core at a time along the chain until the final core is reached. In this process, the output indices are propagated along the chain, while the first core, having lost its outer indices, is absorbed into its neighboring core. This restores the regular TTO form, as illustrated in Fig.~\ref{fig:qft_qc}(c), step~2.

\begin{figure*}
    \centering
    \includegraphics[width=\linewidth]{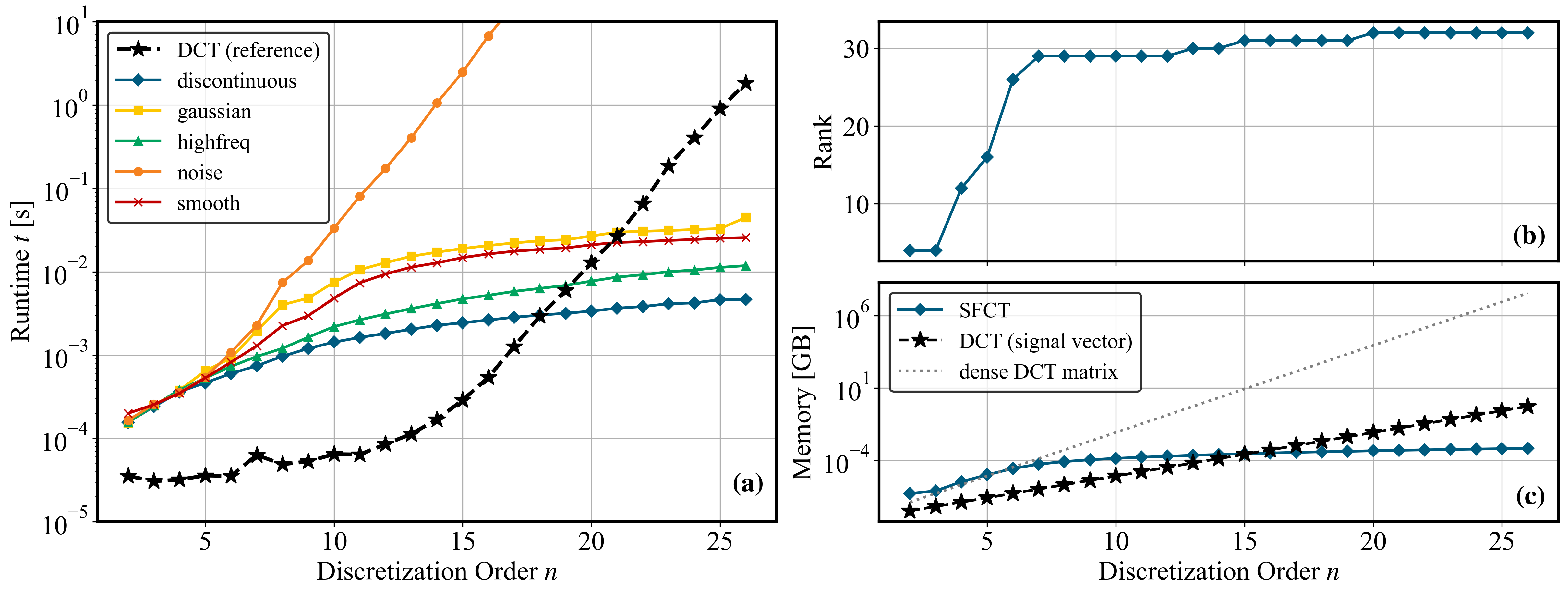}
    \caption{
    (a) Runtime comparison of the classical DCT and the SFCT for different input classes as a function of the discretization order $n$;
    (b) TTO rank evolution of the SFCT;
    (c) Memory of the SFCT operator compared with the signal vector
    required by the classical FFT-based DCT (dashed) and, for reference, the dense DCT matrix (dotted).
    }
    \label{fig:exp_DCTvsSFCT}
\end{figure*}

\subsection{Validation}
\label{subsec:numerical_results}

Before incorporating the superfast spectral transform into practical applications, we first assess the performance of the SFCT by comparing it with the classical DCT for a range of representative input signals. Specifically, we investigate its accuracy, computational cost, and TTO rank growth as functions of the discretization order, characterized by the number of qubits $n$.

To probe different spectral characteristics, we consider the following classes of inputs:
(i) smooth periodic signals,
(ii) high-frequency oscillatory signals,
(iii) discontinuous signals,
(iv) localized Gaussian profiles, and
(v) i.i.d. Gaussian white noise samples.
These test cases span a range of regularity and spectral localization properties, which are known to strongly influence TT compressibility.
The explicit definitions are provided in Appendix~\ref{app:SFCT_Inputs}.
All experiments are conducted on a one-dimensional uniform grids of size $N = 2^n$, with $n$ ranging from $2$ to $26$.
The classical reference solutions are computed using the orthonormal type-II DCT.
However, since the computational cost of the DCT-II does not depend on the specific input, only one representative run is reported.
The TT-based computations employ an error-controlled truncation scheme with tolerance $\varepsilon_\mathrm{TT} = 10^{-8}$.
Consequently, the error between the reference solution and the TT approximation remains on the order of the prescribed tolerance and is therefore not shown in a separate plot. No explicit rank constraints are imposed, allowing the TT ranks to adapt dynamically.

Figure~\ref{fig:exp_DCTvsSFCT}(a) shows the runtime comparison between the classical DCT and the SFCT. While the classical implementation exhibits the expected $\mathcal{O}(N \log N)$ scaling~\cite{Makhoul1980}, the TT-based method shows favorable scaling behavior for all structured input types. In contrast, the unstructured white noise input does not exhibit compression and follows a scaling comparable to the classical method, up to a constant prefactor.
For structured inputs, all TT-based runs exhibit a crossover with the classical method at approximately $n \in [19,21]$. Beyond this regime, the SFCT becomes increasingly efficient, with runtime improvements of several orders of magnitude.

The evolution of the maximal TT rank of the SFCT representation is shown in Figure~\ref{fig:exp_DCTvsSFCT}(b). After an initial growth phase for small $n$, the ranks stabilize at a plateau value of approximately $32$. This behavior is consistent with Theorem~\ref{theo:expressivity_theorem}, under the assumption that the underlying reduced QFT operator saturates at a fixed rank ($r_S = 8$ for the chosen accuracy).

Finally,
Figure~\ref{fig:exp_DCTvsSFCT}(c) compares the memory required by the SFCT
operator with that of the classical DCT. Since the classical
transform is evaluated via the FFT without forming the
transform matrix, its memory is dominated by the signal
vector and grows as $\mathcal{O}(N)=\mathcal{O}(2^n)$, while the dense $2^n\times 2^n$
transform matrix is shown for reference only. Once the TT
rank plateau is reached, the memory of the SFCT grows only
linearly in $n$. The two curves cross between $n=15$ and $n=16$, and at $n=26$ the SFCT operator requires \SI{0.67}{\mega\byte}
compared with \SI{540}{\mega\byte} for the signal vector, a reduction by
a factor of about 800.

Overall, these results demonstrate that the SFCT provides an efficient and accurate representation of spectral transforms for inputs that admit a compact TT representation, particularly in the high-resolution regime of large $n$.
The corresponding results for the DST-II and the SFST are
given in Appendix~\ref{app:exps_dct2}.

\section{Application I: Quantum-Inspired CFD}
\label{sec:quantum_inspired_cfd}

TT methods have been used to simulate incompressible flows in
periodic~\cite{Gourianov2022,Hoelscher2025,Pisoni2026} and wall-bounded
domains~\cite{Kiffner2023,Peddinti2024,VANHULST2026110169}.
Several of these methods solve the Poisson equation for pressure or streamfunction in physical space.
For the two-dimensional Rayleigh--B\'enard benchmark of Ref.~\cite{vanHulst2026RBC}, we compare two spectral TT pressure solvers with a preconditioned physical-space TT solver. In this problem, a fluid layer is heated
from below and cooled from above~\cite{Ahlers2009}.
The horizontal direction is periodic, and at the plates we impose no-slip
conditions on the velocity and homogeneous Neumann conditions in the
pressure projection~\cite{Chorin1979,Guermond2006}. The discrete pressure operator is diagonalized by a Fourier transform in the horizontal direction and the DCT-II in the wall-normal direction (Section~\ref{sec:spectral_methods}).

Following Ref.~\cite{vanHulst2026RBC}, we use
the Oberbeck--Boussinesq approximation and scale lengths by the plate
separation $L$, temperatures by their difference $\Delta$, and velocities
by $U_{\mathrm{ff}}=\sqrt{g\beta\Delta L}$, where $g$ is the gravitational
acceleration and $\beta$ the thermal expansion coefficient.
Time is measured in units of $t_0=L/U_{\mathrm{ff}}$.
The Rayleigh and Prandtl numbers are
$Ra=g\beta\Delta L^3/(\nu\kappa)$ and $Pr=\nu/\kappa$, with
kinematic viscosity $\nu$ and thermal diffusivity $\kappa$.
The governing equations are
\begin{align}
    \partial_t u_i+u_j\partial_j u_i
    &= -\partial_i p+\sqrt{\frac{Pr}{Ra}}\,\partial_j\partial_j u_i
       -\theta\delta_{iy},\\
    \partial_t\theta+u_j\partial_j\theta
    &= \frac{1}{\sqrt{Ra\,Pr}}\,\partial_j\partial_j\theta,\\
    \partial_j u_j &= 0.
\end{align}
Here $u_x=u$ and $u_y=v$, repeated indices $i,j\in\{x,y\}$ are summed,
and $\delta_{iy}$ is the Kronecker delta. The coordinate $y$ increases
downward, as in Ref.~\cite{vanHulst2026RBC}. The pressure $p$ is scaled
by the reference density times $U_{\mathrm{ff}}^2$, and $\theta$ denotes
temperature relative to the plate-mean temperature.
We use $Ra=10^8$ and $Pr=1$ in a domain of width $2L$.
At the plates $y=\pm1/2$ in scaled
coordinates, $u=v=0$ and $\theta=\pm1/2$.
The figures use height above the lower plate, $z/L=1/2-y$, and
upward velocity $w=-v$.

We use the uniform staggered grid, second-order finite differences,
and two-stage implicit--explicit Runge--Kutta (IMEX-RK2) scheme
described in Ref.~\cite{vanHulst2026RBC}. The discretization orders are
$n_x=9$ and $n_y=8$, giving $N_x=2^{n_x}=512$ and
$N_y=2^{n_y}=256$ grid points, and $\Delta t/t_0=1/512$.
Let $D_h$ and $\nabla_h$
denote the staggered-grid divergence and gradient. Each pressure
projection solves
\begin{equation}
    D_h\nabla_h p
    = \frac{D_h\boldsymbol u^*}{\Delta t_{\mathrm p}},
    \label{eq:rb_pressure_poisson}
\end{equation}
where $\boldsymbol u^*$ is the intermediate velocity and
$\Delta t_{\mathrm p}$ is $\Delta t/2$ in the first stage and $\Delta t$
in the second. The velocity correction is
$-\Delta t_{\mathrm p}\nabla_h p$.

For the solver comparison, we prescribe a rank cap $r$ for velocity,
temperature and pressure instead of a compression tolerance as in
Section~\ref{subsec:numerical_results}. We apply the same cap when forming
the pressure right-hand side (RHS). To assess storage efficiency, we use the
ratio of stored TT coefficients to entries in the full-grid array for each field:
\[
    S_{\mathrm{rel}}=\frac{\sum_{k=1}^{n_x+n_y}2r_{k-1}r_k}{N_xN_y}.
\]
The three pressure solvers are as follows, where both spectral solvers use the SFFT in the horizontal direction but compute the wall-normal DCT-II differently:

\begin{figure*}[!t]
    \centering
    \includegraphics[width=.98\textwidth]{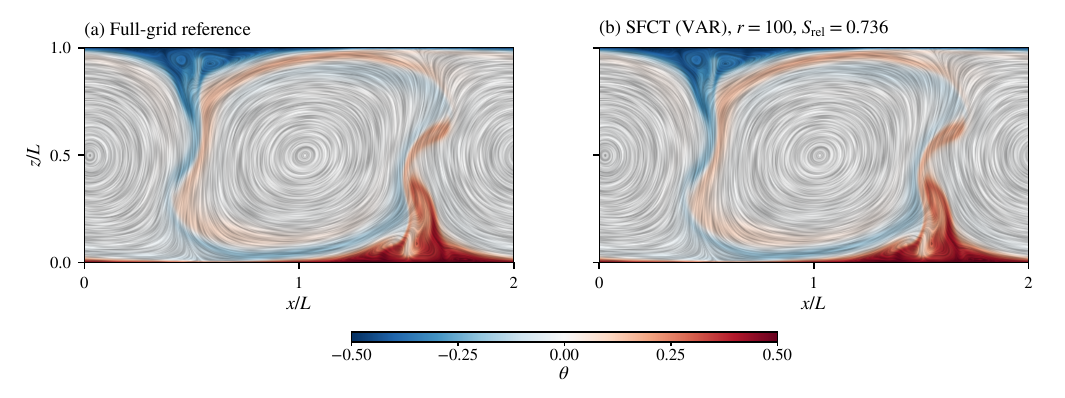}
    \caption{Rayleigh--B\'enard convection at $t/t_0=2$.
    (a) Full-grid reference and (b) SFCT solution at $r=100$, using VAR for the inverse-spectrum product. Parameters are $Ra=10^8$, $Pr=1$, $n_x=9$ and
    $n_y=8$ ($512\times256$ grid points). Color shows temperature and
    line integral convolution (LIC) shows velocity orientation.}
    \label{fig:rb_flow_snapshots}
\end{figure*}

\begin{enumerate}
    \item \emph{Preconditioned ALS.} We solve the physical-space pressure
    equation with ALS~\cite{Holtz2012}, using a left preconditioner
    $P$ with maximum operator rank $r_{\mathrm P}$~\cite{VANHULST2026110169}.
    The construction of $P$ and its use in the local linear systems are
    described in Appendix~\ref{app:preconditioned_als}.
    \item \emph{SFCT.} The SFCT applies the DCT-II
    directly to the $n_y$ wall-normal sites, without extending the field.
    \item \emph{Even extension + SFFT.} The wall-normal DCT-II is
    evaluated through an even extension and the SFFT on $n_y+1$ sites,
    followed by restriction and phase correction. The remaining
    pressure-solve operations are the same as in method 2.
\end{enumerate}
The periodic--Neumann Laplacian has a constant null mode. We subtract
the discrete mean of the RHS before solving, as required by
Eq.~\eqref{eq:poisson_compatibility}.
ALS fixes the pressure at one grid point. The spectral solvers set its
mean to zero. For an exact solve with a compatible RHS, these gauges
give the same pressure gradient. Writing $b$ for the RHS and hats for
the mixed Fourier--cosine transform, the spectral solvers compute
\begin{equation}
    \widehat p=G\odot\widehat b,\qquad
    G_{k\ell}=\begin{cases}
        (\lambda_k^{(x)}+\lambda_\ell^{(y)})^{-1},
            &(k,\ell)\ne(0,0),\\
        0, &(k,\ell)=(0,0).
    \end{cases}
    \label{eq:rb_inverse_spectrum}
\end{equation}
Here $k=0,\ldots,N_x-1$ and $\ell=0,\ldots,N_y-1$, with $(0,0)$
the constant mode. For scaled grid spacings $h_x=2/N_x$ and $h_y=1/N_y$,
the eigenvalues are
\begin{equation*}
 \lambda_k^{(x)}=-\frac{4}{h_x^2}\sin^2\!\frac{\pi k}{N_x},\qquad
 \lambda_\ell^{(y)}=-\frac{4}{h_y^2}\sin^2\!\frac{\pi\ell}{2N_y}.
\end{equation*}
The wall-normal eigenvalues have the cosine spectrum derived in
Section~\ref{sec:spectral_methods}.
For the SFCT, we approximate the inverse-spectrum product using VAR.
For even extension + SFFT, we compare VAR with RSI~\cite{Meng2026-arxiv} for the same product.
We compare these three spectral configurations with preconditioned ALS.
We construct the TT representation of $G$ using tensor cross
interpolation (TCI)~\cite{Oseledets2010B}.
The spectral construction tolerance
$\varepsilon_{\mathrm{spec}}$ controls $G$, the eigenvalue vectors and
phase factors. The inverse transform is
followed by realification (Appendix~\ref{app:realification}) and rounding
to $r$.

Starting from the same developed convective state, we compare
trajectories through $t/t_0=2$.
We run each of the three spectral configurations at $r=50,100,150$ and preconditioned ALS only at $r=100$.
Separate short runs compare the computational cost of all four solver
configurations at $r=100$.

\begin{figure*}[!t]
    \centering
    \includegraphics[width=.98\textwidth]{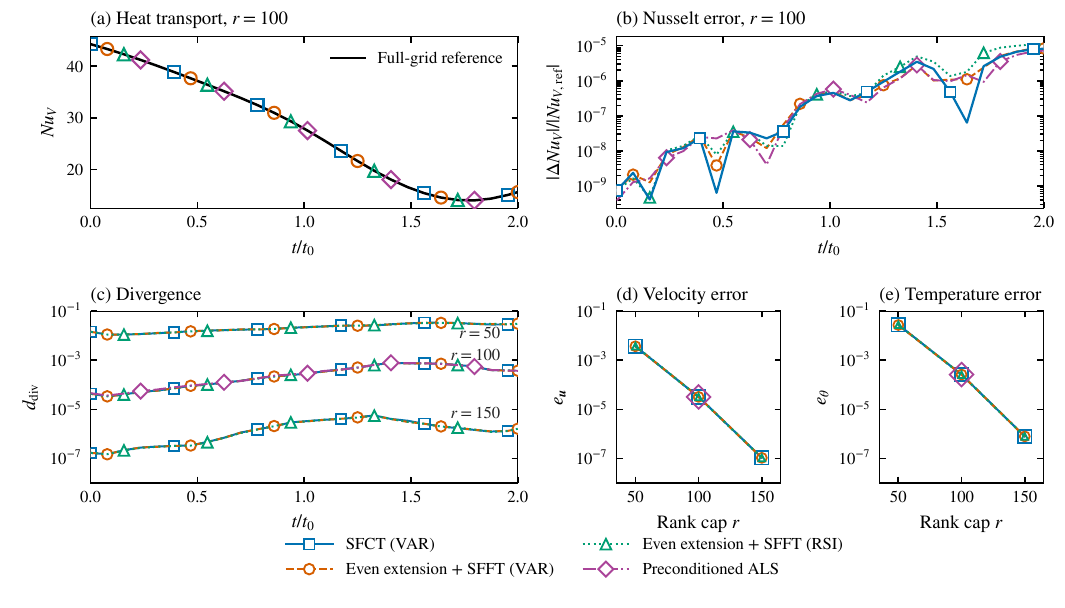}
    \caption{Heat transport, velocity divergence and field errors for the four pressure-solver configurations.
    Spectral constructions use $\varepsilon_{\mathrm{spec}}=10^{-10}$.
    (a) Volume Nusselt number and (b) its relative error at $r=100$.
    (c) Divergence RMS at the indicated ranks.
    (d) Relative velocity error and (e) relative temperature error at $t/t_0=2$.
    Preconditioned ALS is shown only at $r=100$.}
    \label{fig:rb_observables}
\end{figure*}

For velocity and temperature, we use
\begin{equation}
 \begin{aligned}
 e_{\boldsymbol u}&=
 \frac{\|\boldsymbol u-\boldsymbol u_{\mathrm{ref}}\|_2}
      {\|\boldsymbol u_{\mathrm{ref}}\|_2},\\
 e_\theta&=\frac{\|\theta-\theta_{\mathrm{ref}}\|_2}
                    {\|\theta_{\mathrm{ref}}\|_2}.
 \end{aligned}
 \label{eq:rb_field_errors}
\end{equation}
The velocity norm includes both components. The reference uses the
same spatial and temporal discretization without TT compression.
We evaluate heat transport and incompressibility using
\begin{equation}
 \begin{aligned}
 Nu_V(t)&=1-\sqrt{Ra\,Pr}\,\langle v\theta\rangle_V,\\
 d_{\mathrm{div}}(t)&=\frac{\|D_h\boldsymbol u\|_2}{\sqrt{N_xN_y}}.
 \end{aligned}
 \label{eq:rb_observables}
\end{equation}
Here $\langle\cdot\rangle_V$ denotes the volume average.
The minus sign follows from the
downward velocity $v=-w$.

Figure~\ref{fig:rb_flow_snapshots} compares the reference and SFCT
fields at $t/t_0=2$. The rank-100 TT fields have a relative storage
requirement of $S_{\mathrm{rel}}=0.736$.
All four rank-100 runs in Fig.~\ref{fig:rb_observables}(a,b)
reproduce the reference heat transport with
relative deviations below $1.1\times10^{-5}$ over the sampled times.

\begin{figure}[pos=t]
    \centering
    \includegraphics[width=\columnwidth]{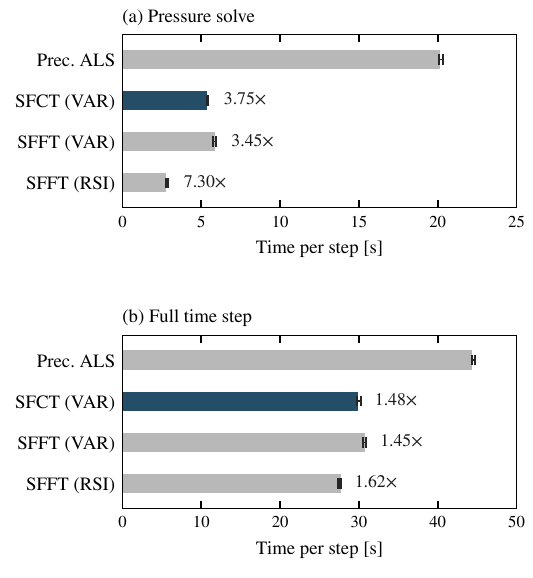}
\caption{Time per step at $r=100$ for (a) both pressure projections combined and (b) the full step. SFFT denotes the even-extension runs. Bars show median times. Error bars span the three repetitions.
 Each label shows the median speedup over preconditioned ALS (Prec.\ ALS), calculated from matched runs.
 See Appendix~\ref{app:rb_comparison} for timing procedures and solver details.}
 \label{fig:rb_matched_accuracy}
\end{figure}

Increasing the field rank reduces both the divergence and the field
errors (Fig.~\ref{fig:rb_observables}(c--e)). The SFCT and even extension followed by the SFFT give similar errors at all three ranks. Their endpoint
divergence RMS falls from about $3.1\times10^{-2}$ at $r=50$ to
$1.6\times10^{-6}$ at $r=150$.
At $r=100$, all four configurations give endpoint velocity and
temperature errors of approximately $3.2\times10^{-5}$ and
$2.6\times10^{-4}$, with divergence RMS values of
\mbox{$(3.74\text{--}3.84)\times10^{-4}$}. At $r=150$, the spectral runs
reduce the field errors to about $1.1\times10^{-7}$ and
$7.9\times10^{-7}$, respectively, roughly 300 times smaller than at $r=100$.
The full-grid reference divergence remains below $8.5\times10^{-15}$.
In the TT runs, divergence can remain even after an accurate pressure solve because we round the RHS and the updated velocity to rank $r$. Changing VAR to RSI
in the SFFT runs with even extension has comparatively little effect on the
field errors.

Figure~\ref{fig:rb_matched_accuracy} compares pressure-solve and full-step times at $r=100$ for the four configurations.
At step 10, all four give relative errors of approximately $7.7\times10^{-7}$ in velocity and $6.2\times10^{-6}$ in temperature.
Relative to preconditioned ALS, the SFCT with VAR gives a $3.75$-fold speedup for the pressure solve and a $1.48$-fold speedup for the full step. For even extension followed by the SFFT, the corresponding speedups are $3.45$ and $1.45$ with VAR, and $7.30$ and $1.62$ with RSI. The momentum and temperature solves account for most of the remaining cost.

\section{Application II: Quantum-Inspired Homogenization}
\label{sec:quantum_inspired_homogenization}
Given a composite material, Fourier-based homogenization~\cite{Moulinec1998, Schneider2021} determines its macroscopic mechanical response by solving a microscopic boundary value problem on a representative volume element (RVE) $\Omega=[0,L]^D$~\cite{Cioranescu1999}. For a prescribed macroscopic strain $\mathbf{E}$, the objective is to determine the effective stiffness tensor $\mathbf{C}_{\mathrm{eff}}$, relating the macroscopic stress and strain through
\begin{align*}
\boldsymbol{\sigma}_{\mathrm{eff}}
= \mathbf{C}_{\mathrm{eff}}:\mathbf{E}.
\end{align*}
The RVE considered here consists of distinct material phases with piecewise-constant properties and sharp interfaces. The local constitutive response is described by
\begin{align}
\boldsymbol{\sigma}(\mathbf{x})
= \mathbf{C}(\mathbf{x}):\boldsymbol{\varepsilon}(\mathbf{x}).
\end{align}
We first review the classical FFT-based homogenization scheme and its Tensor-Train extension before generalizing the latter to periodicity-compatible mixed uniform boundary conditions (PMUBC), which are particularly suitable for finite microstructural samples that are not intrinsically periodic. Their spectral treatment relies on the DCT-II and DST-II formulations introduced above.

Throughout this section, stress and strain tensors are represented in Voigt notation as vectors $\boldsymbol{\sigma},\boldsymbol{\varepsilon}\in\mathbb{R}^{n_v}$, with $n_v=3$ in two dimensions and $n_v=6$ in three dimensions. The stiffness tensor is correspondingly represented by a matrix $\mathbf{C}\in\mathbb{R}^{n_v\times n_v}$~\cite{Schneider2021}.

\subsection{FFT-based Homogenization}

On the periodic RVE $\Omega$, the total displacement field is given by $\mathbf{u}(\mathbf{x})=\mathbf{E}\mathbf{x}+\tilde{\mathbf{u}}(\mathbf{x})$, where $\mathbf{E}\mathbf{x}$ denotes the prescribed macroscopic part and $\tilde{\mathbf{u}}(\mathbf{x})$ is a periodic fluctuation.
The local strain field is obtained from the displacement field as
\begin{align}
    \boldsymbol{\varepsilon}(\mathbf{x}) = \mathbf{E} + \nabla^S \tilde{\mathbf{u}}(\mathbf{x}),
    \label{eq:kinematics}
\end{align}
with $\nabla^S$ denoting the symmetric part of the gradient. We further require the fluctuation to have a vanishing average gradient over the cell,
\begin{align*}
    \frac{1}{|\Omega|}\int_\Omega \nabla\tilde{\mathbf{u}}(\mathbf{x})\,\mathrm{d}\mathbf{x} = \mathbf{0},
\end{align*}
which ensures that $\mathbf{E}$ equals the volume-averaged strain.

Requiring the local stress field to satisfy equilibrium and introducing
a homogeneous reference stiffness $\mathbf{C}_0$ together with the
polarization field
\[
    \tau(\mathbf{x}) =
    [\mathbf{C}(\mathbf{x}) - \mathbf{C}_0]\,\varepsilon(\mathbf{x}),
\]
the local equilibrium problem can be recast as the well-known
Lippmann--Schwinger equation~\cite{Moulinec1998, Brisard2010}
\begin{equation}
    \varepsilon(\mathbf{x}) =
    \mathbf{E} -
    \left\{\Gamma_0 \ast
    \left[(\mathbf{C} - \mathbf{C}_0)\,\varepsilon\right]\right\}(\mathbf{x}).
    \label{eq:lippmann-schwinger}
\end{equation}
The symbol $\Gamma_0$ denotes the Green--Eshelby operator associated
with the homogeneous reference stiffness $\mathbf{C}_0$.
This operator admits a closed-form representation in Fourier space~\cite{Willot2015},
and, additionally, the convolution in the Lippmann--Schwinger equation
reduces to a pointwise product in this representation.
This enables efficient evaluation using the Fast Fourier Transform (FFT),
denoted by $\mathcal{F}$, resulting in the FFT-based homogenization
algorithm~\ref{algo:fft}.
This algorithm evaluates Eq.~\eqref{eq:lippmann-schwinger} iteratively,
starting from $\varepsilon_0 = \mathbf{E}$, until the convergence
criterion $\delta_{\mathrm{tol}}$ is satisfied.
Here and in the following we set $\delta_{\mathrm{tol}}=10^{-4}$.

After convergence, the effective stiffness tensor $\mathbf{C}_{\mathrm{eff}}$
is obtained from the volume-averaged stress corresponding to each
prescribed macroscopic strain $\mathbf{E}$. Specifically,
\begin{align}
    \boldsymbol{\sigma}_{\mathrm{eff}}
    =
    \frac{1}{|\Omega|}
    \int_\Omega
    \mathbf{C}(\mathbf{x})\,\varepsilon(\mathbf{x})\,\mathrm{d}\mathbf{x},
\end{align}
and the resulting stress--strain relation determines
$\mathbf{C}_{\mathrm{eff}}$~\cite{Hill1965}.

The computational cost of the FFT-based scheme scales as $\mathcal{O}(N^D\log N)$ for an $N^D$ discretization. This motivates the Tensor-Train formulation introduced in the following subsection, which aims to reduce the computational and memory requirements for high-resolution RVEs.

\begin{algorithm}[t]
    \SetAlgoLined
    \DontPrintSemicolon

	\KwIn{Macroscopic strain $\mathbf{E}$, stiffness $\mathbf{C}$, tolerance $\delta_{\mathrm{tol}}$}
	\KwOut{Local strain field $\boldsymbol{\varepsilon}$}

    $\boldsymbol{\varepsilon}^{(0)}
    \gets
    \mathbf{E}$\;
    $m \gets 0$\;

    \While{not converged}{
        $\boldsymbol{\tau}^{(m)}
        \gets
        (\mathbf{C}-\mathbf{C}^0)
        \boldsymbol{\varepsilon}^{(m)}$\;

        $\hat{\boldsymbol{\tau}}^{(m)}
        \gets
        \mathcal{F}
        \left[
        \boldsymbol{\tau}^{(m)}
        \right]$\;

        $\hat{\boldsymbol{\varepsilon}}^{(m)}
        \gets
        \hat{\varGamma}^0 
        \,
        \hat{\boldsymbol{\tau}}^{(m)}$\;

        $\boldsymbol{\varepsilon}^{(m+1)}
        \gets
        \mathbf{E}
        -
        \mathcal{F}^{-1}
        \left[
        \hat{\boldsymbol{\varepsilon}}^{(m)}
        \right]$\;

		\BlankLine
		\If{$\displaystyle \left\|\boldsymbol{\varepsilon}^{(m+1)}-\boldsymbol{\varepsilon}^{(m)}\right\|_F<\delta_{\mathrm{tol}} \left\|\mathbf{E}\right\|_F$}{
   		\textbf{break}\;
	}
        $m \gets m+1$\;
    }

    \Return{$\boldsymbol{\varepsilon}^{(m)}$}\;

    \caption{FFT-based homogenization algorithm}
    \label{algo:fft}
\end{algorithm}

\subsection{Tensor-Train based extension}
The FFT-based homogenization scheme can be transferred to the tensor train format by representing the discretized fields and operators as low-rank tensor representations and replacing the classical Fourier transform by the Superfast Fourier Transform (SFFT).
A detailed derivation of the resulting SFFT-based homogenization scheme is given in~\cite{Hauck2026};
here, we summarize the aspects relevant for the present formulation.
All fields and operators use the QTT representation of Section~\ref{sec:tensor_trains}.

In contrast to the scalar fields often considered in conventional
applications of TT formalism, the quantities appearing in the homogenization algorithm carry additional tensorial indices of dimensionality $n_\nu$ associated with the Voigt representation.
The TT representation must therefore account for these indices in addition to the indices associated with the discretization procedure.
This is achieved by introducing an additional \emph{control core} at the beginning of the TT representation, whose dimensionality corresponds to the number of Voigt components, six in the present three-dimensional formulation.
The Green--Eshelby operator requires particular care because its momentum-space representation imposes a specific TT core ordering that must be consistent with the convention dictated by the SFFT (cf.~Section~\ref{sec:spectral_methods}).
For a multidimensional tensor-product transform, this is achieved by reversing the ordering of the TT cores independently within each spatial dimension.

With these representations, the FFT-based homogenization procedure can
be transferred directly to TT format by replacing the Fourier
transforms $\mathcal{F}$ and $\mathcal{F}^{-1}$ with the SFFT and its
inverse, denoted by $\mathcal{S}$ and $\mathcal{S}^{-1}$, respectively.
The pointwise application of the Green--Eshelby operator,
$\hat{\Gamma}_0\hat{\tau}_m$, is replaced by the fast AMEn TT-TTO contraction algorithm denoted by $\mathcal{A}(\hat{\Gamma}_0,\hat{\tau}_m)$~\cite{Dolgov2013}.

The resulting SFFT-based homogenization scheme therefore follows the
same fixed-point iteration as the classical FFT formulation.
After convergence, the local strain field is obtained
by contracting the resulting tensor train representation.
If the TT ranks of all fields and operators remain bounded, the computational cost per iteration grows only polylogarithmically with $N$, in contrast to the $\mathcal{O}(N^D \log N)$ complexity of the classical FFT-based scheme~\cite{Hauck2026}.
Whether this regime is reached depends on the microstructure.

\subsection{Spectral-Tensor-Train Homogenization}
\label{sec:stt-homogenization}

We generalize the SFFT-based homogenization framework to Spectral-Tensor-Train (STT) homogenization by replacing the purely periodic boundary conditions with periodicity-compatible mixed uniform boundary conditions (PMUBC)~\cite{Wiegmann1999, Grimm2021}.
These conditions extend periodic boundary conditions by combining displacement and traction constraints and are particularly suitable for finite microstructural samples that are not intrinsically periodic.
In particular, for porous media, finite-volume representations of the microstructure may exhibit pronounced size and boundary effects, and artificially periodic continuation can alter features such as pore connectivity and therefore lead to significant differences in the predicted effective properties~\cite{Gelebart2020, Pahr2008, Risthaus2024dirichlet}.
For readability, we present the following construction in two dimensions. Its extension to three dimensions is straightforward and is adopted in the numerical experiments.

For a prescribed macroscopic tensile load in voigt notation $\mathbf{E}=E_1\mathbf{e}_1$ ,
the corresponding PMUBC are given by
\begin{align*}
    u_1(0,y) &= -\frac{E_1 L_1}{2},
    & u_1(L_1,y) &= \frac{E_1 L_1}{2},
    \\
    \sigma_{3}(0,y) &= 0,
    & \sigma_{3}(L_1,y) &= 0,
    \\
    \sigma_{3}(x,0) &= 0,
    & \sigma_{3}(x,L_2) &= 0,
    \\
    u_2(x,0) &= 0,
    & u_2(x,L_2) &= 0.
\end{align*}
Here, the normal displacement is prescribed in the loading direction,
while the transverse displacement and shear tractions are constrained
according to the PMUBC prescription.

For a macroscopic shear load
$\mathbf{E}=E_3\mathbf{e}_3$, the corresponding conditions take a different form, fulfilling
\begin{align*}
    \sigma_{1}(0,y) &= 0,
    & \sigma_{1}(L_1,y) &= 0,
    \\
    u_2(0,y) &= -\frac{E_3 L_1}{4},
    & u_2(L_1,y) &= \frac{E_3 L_1}{4},
    \\
    u_1(x,0) &= -\frac{E_3 L_2}{4},
    & u_1(x,L_2) &= \frac{E_3 L_2}{4},
    \\
    \sigma_{2}(x,0) &= 0,
    & \sigma_{2}(x,L_2) &= 0.
\end{align*}
Here, the roles of displacement and traction constraints are exchanged
component-wise. The particular choice of boundary condition therefore
depends on the applied macroscopic strain component, rather than being
identical for all load cases. Details for all components in the
three-dimensional formulation are given in
Appendix~\ref{app:PMUBC}.

PMUBC can alternatively be interpreted through an extended mirrored
domain. In particular, PMUBC on the original domain $\Omega=[0,L]^2$ are
equivalent to periodic boundary conditions on
$\widetilde{\Omega}=[-L,L]^2$, provided that the
displacement field is extended according to
\begin{align*}
    u_1(-x,y) &= -u_1(x,y),
    & u_2(-x,y) &= u_2(x,y),
    \\
    u_1(x,-y) &= u_1(x,y),
    & u_2(x,-y) &= -u_2(x,y).
\end{align*}
These reflection symmetries induce corresponding even and odd symmetries in the local strain and stress components.
The resulting extensions satisfy the appropriate matching conditions across the boundaries of the mirrored cells and therefore define periodic fields on the extended domain~\cite{Wiegmann1999}.

\begin{figure*}
    \centering
    \includegraphics[width=0.75\linewidth]{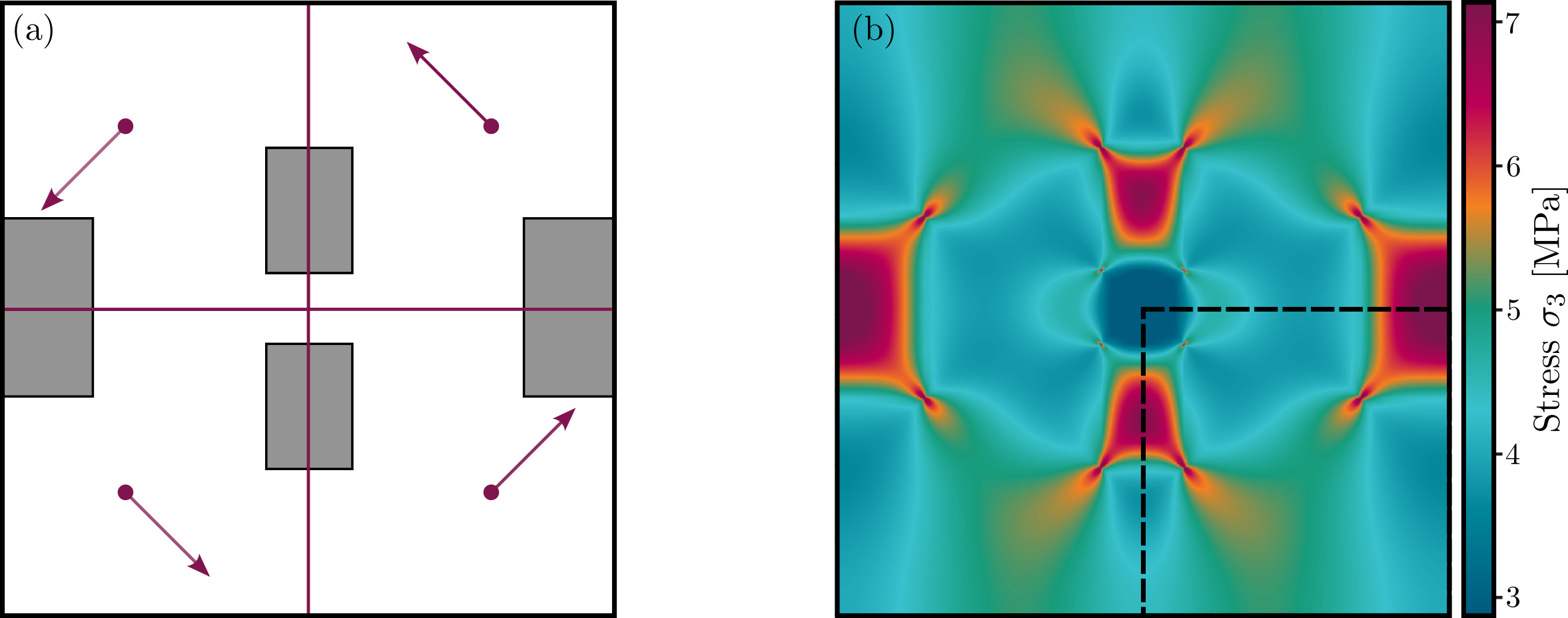}
    \caption{Slice of the 3D RVE (a) for the input geometry and (b) the corresponding local stress field $\sigma_3$ under shear load.
    The arrows indicate the relation between the displacement vectors within each quadrant.
    The dashed box denotes the reduced domain used for the STT homogenization method.
    The gray region has a Young's modulus of $E_1=\SI{20}{\giga\pascal}$,
    while the white region has $E_2=\SI{1}{\giga\pascal}$.
    Both regions have a Poisson's ratio of $\nu_1=\nu_2=1/3$.}
    \label{fig:mirrored_domain}
\end{figure*}

Although this construction permits the direct application of the
FFT-based homogenization scheme on the mirrored domain, it increases
the number of grid points by a factor of four in two dimensions and
eight in three dimensions. To avoid this severe overhead, we instead exploit
the symmetries induced by the mirrored-domain construction directly on
the original domain, following~\cite{Grimm2021}. In particular, each field
component is either even or odd with respect to a given coordinate,
with the specific symmetry depending on the applied macroscopic strain
and the field component. This structure allows the corresponding
Fourier transforms to be evaluated using only half of the original
grid.

The construction parallels the even extension of Section~\ref{sec:spectral_methods}, with two differences: the symmetry is now imposed
by the boundary conditions rather than chosen, and it differs
between field components and load cases.
Consider first a discretized volumetric component
$f_j$, $j=0,\ldots,N-1$, for which the even or odd symmetry on the
mirrored domain translates to the grid-level relation
\begin{align*}
    f_j = \pm f_{N-1-j},
    \qquad
    \begin{cases}
        + & \text{even symmetry},\\
        - & \text{odd symmetry}.
    \end{cases}
\end{align*}
Substituting this relation into the formula for the discrete Fourier coefficients
\begin{align*}
    \hat{f}_k =
    \sum_{j=0}^{N-1}
    f_j \exp\left(\frac{2\pi \mathrm{i}jk}{N}\right)
\end{align*}
and pairing each index $j<N/2$ with its mirrored counterpart
$N-1-j$ yields, for an even number of grid points,
\begin{align*}
    \hat{f}_k^{\mathrm{even}}
    &=
    2\exp\left(-\frac{\mathrm{i}\pi k}{N}\right)
    \sum_{j=0}^{N/2-1}
    f_j
    \cos\left(\frac{\pi(2j+1)k}{N}\right),
    \\
    \hat{f}_k^{\mathrm{odd}}
    &=
    2\mathrm{i}\exp\left(-\frac{\mathrm{i}\pi k}{N}\right)
    \sum_{j=0}^{N/2-1}
    f_j
    \sin\left(\frac{\pi(2j+1)k}{N}\right).
\end{align*}
The phase factor accounts for the fact that the reflection symmetry is
centered between grid points, while the remaining sums are real-valued
cosine and sine transforms~\cite{Schwarztrauber1986, Grimm2021}.

For shear components, the corresponding symmetry relation takes the
form $f_j=\pm f_{N-j}$, resulting in
\begin{align*}
    \hat{f}_k^{\mathrm{even}}
    &=
    f_0 + (-1)^k f_{N/2}
    + 2\sum_{j=1}^{N/2-1}
    f_j\cos\left(\frac{2\pi jk}{N}\right),
    \\
    \hat{f}_k^{\mathrm{odd}}
    &=
    2\mathrm{i}\sum_{j=1}^{N/2-1}
    f_j\sin\left(\frac{2\pi jk}{N}\right).
\end{align*}
Thus, the PMUBC formulation requires discrete cosine and sine
transforms with additional phase factors, which are realized
in the TT format by the SFCT and SFST, respectively.

Depending on the applied macroscopic strain, different components of
the displacement, strain, and stress fields exhibit different symmetry
types and consequently require different combinations of SFCT and SFST, as indicated in Figure~\ref{fig:mirrored_domain}.
The complete component-wise transformation pattern for the employed
Hex8R discretization~\cite{Schneider2017} therefore relies exclusively on the
DCT-II and DST-II and is given in Appendix~\ref{app:PMUBC}.
Consistent with the type-II transforms, whose symmetry axis lies between two samples, the microstructure is sampled at the voxel centers, so that the mirror planes coincide with the voxel boundaries.

When assembling these component-wise transforms in the TT format,
the different frequency sets of the cosine and sine transforms require
a common frequency layout.
However, for a transform length $M=N/2$, the cosine
transform yields the frequencies $k=0,\ldots,M-1$, whereas the sine
transform yields $k=1,\ldots,M$. We therefore store each sine coefficient
in slot $k\bmod M$, mapping the Nyquist frequency $k=M$ to slot $0$.
This cyclic shift introduces no information loss, since the sine
coefficient at $k=0$ vanishes identically. However, slot $0$ then
represents different physical frequencies for cosine and sine components.
This distinction must be retained when applying the Green--Eshelby
operator, which couples the Voigt components: at slot $0$, cosine
components are coupled at $k=0$ and sine components at $k=M$, whereas
components of different parity are decoupled. All other frequency
slots remain unchanged.

Returning to Algorithm~\ref{algo:fft}, PMUBC can be incorporated into the Fourier-based homogenization algorithm by replacing the standard SFFT with the load-dependent, component-wise SFCT or SFST and the associated phase factors as well as the periodic Green operator by the modified operator described above. The multidimensional transforms are constructed separably by applying the prescribed one-dimensional transforms along each spatial direction for every Voigt component and load case, according to the symmetry patterns listed in Appendix~\ref{app:PMUBC}. The corresponding one-dimensional TTOs are concatenated along the consecutive directional TT cores and connected by rank-one bonds.

The phase factors are also represented in TT format with rank one. They depend only on the mode index $k$, defined analogously to Eq.~\eqref{eq:binary_index}, along the corresponding direction.
The \(l\)-th core of the phase TT satisfies \(v_l(i_l=0)=1\), while for \(i_l=1\) it is given by
\begin{align*}
v_l(i_l=1)=\exp\left(-\mathrm{i}\frac{\pi\,2^l}{2N}\right).
\end{align*}
The resulting rank-one phase representation can be lifted to a TTO through Eq.~\eqref{eq:tt_to_tto} and combined directly with the corresponding SFCT or SFST TTO.

\subsection{Evaluation of the Homogenization Algorithm}
\label{subsec:numerical_evaluation}

\begin{figure*}[t]
    \centering
    \includegraphics[width=\linewidth]{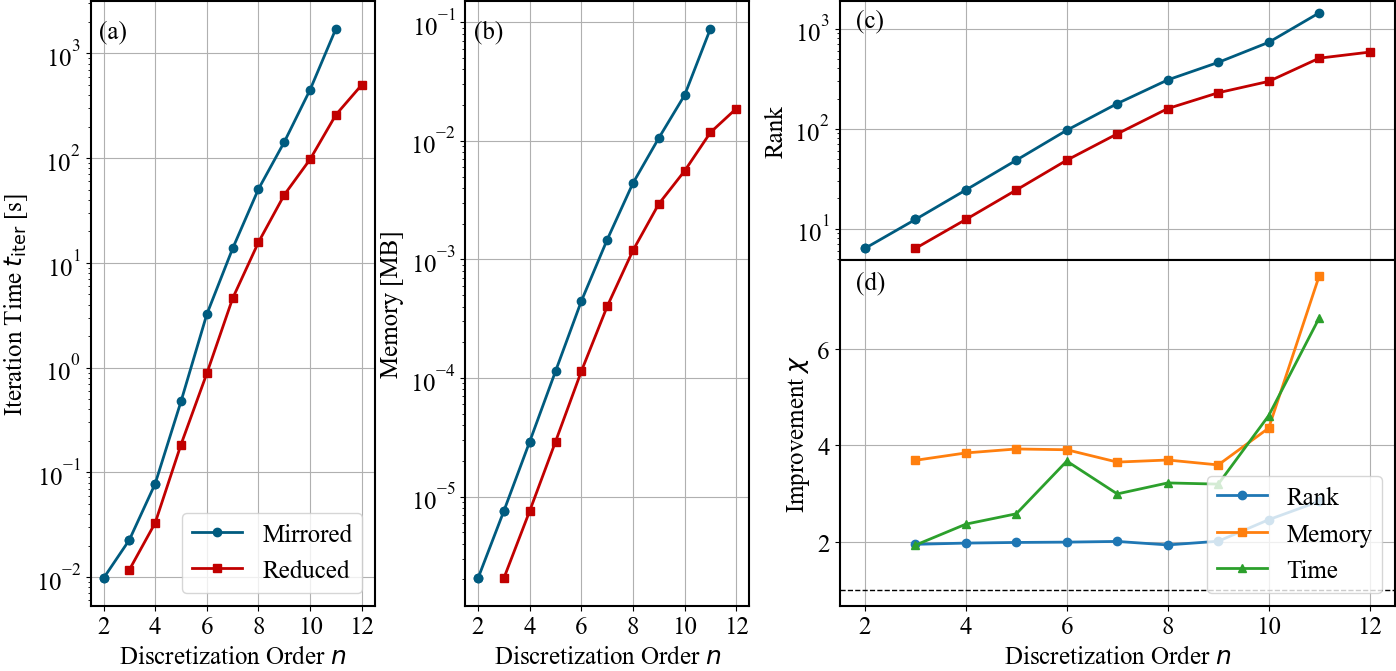}
    \caption{Computational scaling of the mirrored and reduced-domain
    approaches for the three-dimensional boxed-fiber geometry:
    (a) iteration time, (b) memory consumption, (c) maximum TT rank, and
    (d) improvement factor $\chi$ of the reduced approach relative to the
    mirrored approach.}
    \label{fig:homogenization_results}
\end{figure*}

Having established the role of the SFCT and SFST within the STT-based
homogenization algorithm, we evaluate the computational performance of the
mirrored and reduced-domain formulations for the three-dimensional
boxed-fiber geometry shown in Figure~\ref{fig:mirrored_domain}. The
discretization order $n$ determines $N=2^n$ discretization points per
spatial dimension and $N^3 = 2^{3n}$
points in the full three-dimensional domain. Thus, the maximum investigated
order of $n=12$ corresponds to $N=4096$ points per dimension and
approximately $6.87\times10^{10}$ points in the full domain. For the
reduced-domain formulation, the smaller computational domain is explicitly
accounted for in the PMUBC formulation, where a mirrored domain of order $n$ is represented by a reduced domain of order $n-1$, such that both approaches are evaluated at the same physical resolution. All
experiments use a fixed TT-compression tolerance of
$\varepsilon_\mathrm{TT}=10^{-8}$.
The TT-based homogenization reproduces the effective stiffness
of the classical FFT-based scheme over the range in which the latter remains feasible, see Appendix~\ref{app:validation}.

The reduced-domain formulation consistently lowers the computational cost
over the investigated discretization range,
as can be seen in Figure~\ref{fig:homogenization_results}.
To quantify the reduction, we define the improvement factor as
\begin{equation}
    \chi = \frac{Q_{\mathrm{mirrored}}}{Q_{\mathrm{reduced}}},
    \label{eq:improvement_factor}
\end{equation}
where $Q$ denotes memory consumption, iteration time, or maximum TT rank.
Thus, $\chi>1$ indicates an advantage of the reduced formulation. The
memory consumption is reduced by approximately a factor of $4$ over a
substantial range of discretization orders, while the iteration-time
improvement ranges from approximately $2$ to $4$ for $3\leq n\leq10$, matching the range previously reported for the dense setting~\cite{Grimm2021}.
For the highest discretizations, both improvements increase further,
reaching approximately $\chi\approx7$, see Figure~\ref{fig:homogenization_results}(d).
The maximum TT rank is approximately halved by the reduced formulation,
consistent with the observed reduction in memory consumption.

The number of iterations required for convergence remains approximately
the same for both formulations, with only minor variations at the lower
discretization orders. Hence, the reduction in iteration time primarily
results from the lower computational cost per iteration rather than from
a reduction in the number of iterations. Overall, the reduced-domain
formulation provides a consistent computational advantage while maintaining
the same physical discretization resolution.

\section{Conclusion}
\label{sec:conclusion}

We have developed tensor train (TT) representations of the type-II discrete cosine and sine transforms (SFCT and SFST). To this end, we converted a quantum circuit for the combined quantum cosine-sine transform into a single tensor train operator from which both transforms can be extracted. This construction also enabled us to derive an expressivity theorem establishing an upper bound of \(6r_S\) on the TT ranks of the individual transforms, independently of the discretization order. Here, \(r_S\) denotes the TT rank of the Superfast Fourier Transform (SFFT), which saturates at approximately \(r_S \approx 10\). Numerical experiments further show that the ranks of both the SFCT and SFST saturate at approximately \(32\). For structured inputs, both transforms exhibit favorable runtime and memory scaling, outperforming their classical counterparts in runtime for discretization sizes above \(2^{20}\)  (approx. 1 million) and achieving improved memory efficiency above \(2^{15}\) (approx. 32 thousand)  discretization points.

For Rayleigh--B\'enard convection at $r=100$, all four pressure-solver configurations gave comparable velocity and temperature errors. The spectral trajectories also closely reproduced the reference Nusselt number. With VAR for the Hadamard product, the SFCT and even extension followed by the SFFT sped up the Poisson solve by factors of $3.75$ and $3.45$, respectively, relative to preconditioned ALS. Using RSI instead of VAR in the even-extension solver gave a $7.3$-fold speedup for the Poisson solve and a $1.6$-fold speedup for the full time step.
The momentum and temperature solves now account for most of the time per step. We therefore plan to test whether applying recursive sketched interpolation
(RSI)~\cite{Meng2026-arxiv} and alternating cross interpolation
(ACI)~\cite{Ritter2026ACI} to the nonlinear products in the momentum
and temperature equations can further reduce the full-step time.

In computational homogenization, the SFCT and SFST enabled the direct implementation of periodicity-compatible uniform boundary conditions (PMUBC) in TT format for high-resolution three-dimensional problems with up to 70 billion grid points.
Compared with the conventional SFFT-based homogenization algorithm, this direct implementation resulted in a 2--6-fold speedup and a 4--7-fold reduction in memory consumption.

Several extensions of the proposed approach are possible. Other types of discrete cosine and sine transforms can be treated analogously. Moreover, since the Chebyshev transform on Chebyshev--Gauss nodes is equivalent to a DCT-II, the SFCT can also be applied to Chebyshev spectral methods on nonuniform grids. The quantum-circuit-based construction additionally provides a classical framework for simulating and benchmarking quantum spectral algorithms with non-periodic boundary conditions~\cite{Febrianto2026}. More broadly, the SFCT and SFST extend the range of spectral transforms that can be efficiently represented in TT format, opening up new possibilities for high-resolution PDE solvers with non-periodic boundary conditions.

\section{Acknowledgment}
S.H.H., M.K. and N.R.G. gratefully
acknowledge funding from the German National High Performance
Computing (NHR) association for the Center NHR South-West.
N.-L.v.H. and D.J. acknowledge support from the European Union's Horizon Europe
programme under grant agreement No.~101080085 (QCFD).
Their work on the spectral pressure solvers was carried out as part of the 16th German-Israeli
Call for Proposals for Joint R\&D Projects -- Eureka Network.
N.-L.v.H. and D.J. acknowledge computing resources provided by the Hummel-2
cluster at the University of Hamburg, funded by the Deutsche Forschungsgemeinschaft
(DFG, German Research Foundation) -- 498394658.
N.-L.v.H. also thanks Tomohiro Hashizume for early discussions that motivated exploring discrete cosine transforms for the pressure projection.

\appendix
\section{Proof of the Expressivity Theorem}
\label{app:proof_expressivity_theorem}

Let $U$ act on $\mathcal{H}_A\otimes\mathcal{H}_B$ with
$\dim(\mathcal{H}_A\otimes\mathcal{H}_B)=2^n$. For the bipartition between
qubits $j$ and $j+1$, the operator Schmidt decomposition reads
\begin{align}
    U = \sqrt{2^n}\sum_{k=1}^{\hat r}\sigma_k\,A_k\otimes B_k ,
    \label{eq:OSD}
\end{align}
with
$\operatorname{Tr}\big(A_k^{\dagger}A_{k'}\big)=
\operatorname{Tr}\big(B_k^{\dagger}B_{k'}\big)=\delta_{k,k'}$
and $\sum_k\sigma_k^2=1$. The number $\hat r$ of non-zero coefficients is the
operator Schmidt rank; it equals the TT rank across that bipartition, so the
maximum over all bipartitions is the bond dimension of an exact TT
representation. All operators below refer to one such fixed bipartition, and
its index is suppressed; the cases in $j$ are stated explicitly where the rank
depends on it.

Two bounds are used throughout. An expansion into $m$ product operators
$U=\sum_{k=1}^{m}A_k\otimes B_k$ gives
\begin{align}
    \operatorname{rank}_{S}(U)\leq m ,
    \label{eq:rank_product_bound}
\end{align}
without requiring orthogonality of the factors, and the rank is
submultiplicative,
\begin{align}
    \operatorname{rank}_{S}(UV)\leq\operatorname{rank}_{S}(U)\operatorname{rank}_{S}(V).
    \label{eq:rank_submultiplicative}
\end{align}
The only external input is the transform itself: Chen et al. \cite{Chen2023}
show that the Schmidt coefficients of the reduced QFT $\mathcal{Q}_n$ obey
\begin{align*}
    \lim_{n\to\infty}\sigma_k(\mathcal{Q}_n)
    \leq
    \frac{1}{\sqrt{k}}
    \exp\left[-\tfrac{2k+1}{2}\log\left(\tfrac{4k+4}{e\pi}\right)\right],
\end{align*}
so $\mathcal{Q}_n$ admits the low-rank TT representation known as the SFFT,
whose rank we denote by $r_S$. It remains to bound the blocks to the left and
right of the QFT in Figure~\ref{fig:qft_qc}(b).

\subsection{Schmidt ranks of the individual circuit components}
\label{app:individual_components}

Throughout, the controlled permutation following the QFT is the bit-reversed
shift $\widetilde{P}_n^{-1}$ of
Section~\ref{subsec:quantum_spectral_methods}, the convention induced
downstream of the reduced transform $\mathcal{Q}_n$. It enters via the
recursion
\begin{equation}
\begin{aligned}
\widetilde{P}_{n}^{-1} =
    &\Big[\ket{0}\bra{0}^{\otimes j}\otimes\mathds{1}^{\otimes(n-j)}
    +\ket{1}\bra{1}^{\otimes j}\otimes \widetilde{P}_{n-j}^{-1}\Big]\\
    &\times\left(\widetilde{P}_j^{-1}\otimes\mathds{1}^{\otimes(n-j)}\right),
    \end{aligned}
    \label{eq:perm_recursion}
\end{equation}
in which the leading $j$ qubits carry the low-order bits, so that a carry
either stays within the first $j$ qubits or propagates past the cut.

The block to the left of the QFT consists of the Hadamard gate and the
multi-target CNOT,
\begin{align*}
\begin{aligned}
U_1 =
\Big(
 &\ket{0}\bra{0}\otimes\mathds{1}^{\otimes(n-1)}
 +\ket{1}\bra{1}\otimes X^{\otimes(n-1)}
\Big)\\
&\times\left(H\otimes\mathds{1}^{\otimes(n-1)}\right),
\end{aligned}
\end{align*}
which is already of the form \eqref{eq:OSD} across the cut after qubit $j$,
with
\begin{align*}
A_1 &= \tfrac{1}{\sqrt{2^{\,j-1}}}\big(\ket{0}\bra{0}H\big)\otimes\mathds{1}^{\otimes(j-1)},\\
A_2 &= \tfrac{1}{\sqrt{2^{\,j-1}}}\big(\ket{1}\bra{1}H\big)\otimes X^{\otimes(j-1)},\\
B_1 &= \tfrac{1}{\sqrt{2^{\,n-j}}}\,\mathds{1}^{\otimes(n-j)},\\
B_2 &= \tfrac{1}{\sqrt{2^{\,n-j}}}\,X^{\otimes(n-j)} .
\end{align*}
Both families are orthonormal, since $\ket{0}\bra{0}\cdot\ket{1}\bra{1}=0$ and
$\operatorname{Tr}(X^{\otimes(n-j)})=0$, and the two coefficients are equal,
$\sigma_1=\sigma_2=1/\sqrt2$. The left-hand block therefore has
$\operatorname{rank}_{S}(U_1)=2$ across every bipartition.

Right of the QFT, the $K/L/C$ block $U_2$ collapses further. Because
$K_l=\overline{\omega}^{\,2^{l-1}}L_l$, its two controlled branches differ only
by a global phase, so the controlled structure disappears and $U_2$ becomes a
single product operator. Collecting the action on the last qubit into
\begin{align*}
    \bar{C}:=\big(\ket{0}\bra{0}+\overline{\omega}^{\,2^{n-1}-1}\ket{1}\bra{1}\big)C ,
\end{align*}
the cut after qubit $j$ gives the single pair
\begin{align*}
A_1 &= \tfrac{1}{\sqrt{2^{\,j}}}\;L_{1}\otimes\dots\otimes L_{j},\\
B_1 &= \tfrac{1}{\sqrt{2^{\,n-j}}}\;
 L_{j+1}\otimes\dots\otimes L_{n-1}\otimes\bar{C} ,
\end{align*}
with $\sigma_1=1$ and hence $\operatorname{rank}_{S}(U_2)=1$.

The two controlled permutations are identical,
\begin{align*}
    U_3=U_5=\mathds{1}^{\otimes(n-1)}\otimes\ket{0}\bra{0}
    +\widetilde{P}_{n-1}^{-1}\otimes\ket{1}\bra{1},
\end{align*}
and inserting \eqref{eq:perm_recursion} separates them after qubit $j$ into
\begin{equation}
\begin{aligned}
A_1 &= \ket{1}\bra{1}^{\otimes j}\widetilde{P}_j^{-1},
&
B_1 &= \widetilde{P}_{n-j-1}^{-1}\otimes\ket{1}\bra{1},\\
A_2 &= \ket{0}\bra{0}^{\otimes j}\widetilde{P}_j^{-1},
&
B_2 &= \mathds{1}^{\otimes(n-j-1)}\otimes\ket{1}\bra{1},\\
A_3 &= \mathds{1}^{\otimes j},
&
B_3 &= \mathds{1}^{\otimes(n-j-1)}\otimes\ket{0}\bra{0}.
\end{aligned}
\label{eq:U3_expansion}
\end{equation}
The three terms correspond to a carry crossing the cut, a carry absorbed
before it, and the uncontrolled branch. Their left factors are mutually
orthogonal and the terms are independent, so
$\operatorname{rank}_{S}(U_3)=\operatorname{rank}_{S}(U_5)=3$ for
$1\leq j\leq n-2$, reducing to two at $j=n-1$, where the recursion terminates.
The associated Schmidt coefficients are not all equal, so
\eqref{eq:U3_expansion} fixes the rank without being normalised in the sense of
\eqref{eq:OSD}.

The remaining block $U_4$ combines a multi-target CNOT, the gate $B$ on the
last qubit and an inverted multi-control $J$. Writing $D:=(J-\mathds{1})B$ and
$P:=\ket{0}\bra{0}^{\otimes(n-1)}$, the multi-control structure collapses to
$\mathds{1}^{\otimes(n-1)}\otimes B+P\otimes D$, and multiplying by the leading
CNOT gives the four-term expansion
\begin{align*}
U_4 =\;
&\mathds{1}^{\otimes(n-1)}\otimes B\ket{0}\bra{0}
+P\otimes D\ket{0}\bra{0}\\
+\;&X^{\otimes(n-1)}\otimes B\ket{1}\bra{1}\\
+\;&P\,X^{\otimes(n-1)}\otimes D\ket{1}\bra{1} .
\end{align*}
Its four left factors remain linearly independent for $j\le n-2$, so the bound
\eqref{eq:rank_product_bound} is attained and
$\operatorname{rank}_{S}(U_4)=4$, again reducing to two at $j=n-1$. The
explicit Schmidt pairs of $U_4$ are not needed, since only its rank enters
below and the exact rank of the assembled right-hand block is obtained
directly.

\subsection{Rank of the complete block to the right of the QFT}
\label{app:right_block}
Applying \eqref{eq:rank_submultiplicative} to the four components subcircuit $V:=U_2U_3U_4U_5$ yields only
the loose estimate
\begin{align}
    \operatorname{rank}_{S}(V)\leq1\cdot3\cdot4\cdot3=36 ,
    \label{eq:loose_estimate}
\end{align}
which ignores that the components share operator factors: both permutations
are the same $S:=\widetilde{P}_{n-1}^{-1}$, and the projector $P$ appears in
$U_4$ as well. We therefore decompose $V$ directly.

Being a product operator, $U_2$ does not change the rank and is omitted. The
last qubit controls both permutations and carries $B$ and $J$, so it is
natural to expand in its computational basis. Writing
$b_{ab}:=\bra{a}B\ket{b}$ and $d_{ab}:=\bra{a}D\ket{b}$ for the resulting
scalars, one obtains $V=\sum_{a,b}W_{ab}\otimes\ket{a}\bra{b}$ with
\begin{align*}
    W_{ab}=S^a\left(b_{ab}T_b+d_{ab}PT_b\right)S^b
\end{align*}
with the $T_b$ matrices defined through $T_0=\mathds{1}^{\otimes(n-1)}$ and $T_1=X^{\otimes(n-1)} .$
The exponents $a$ and $b$ record whether the permutation acts on the outgoing
and incoming branch, so each $W_{ab}$ is a shift, a gate, and a shift again.

Across the cut after qubit $j$ the factors that act as tensor products split
trivially, $\mathds{1}=\mathds{1}_A\otimes\mathds{1}_R$, $P=P_A\otimes P_R$
and $X^{\otimes(n-1)}=X_A\otimes X_R$. The shift does not: a carry may cross
the cut, which happens exactly when the first $j$ qubits are all one and
become all zero. This gives the two-term split
\begin{align*}
    S=S_A\otimes\mathds{1}_R+E_A\otimes\bar{S}_R ,
\end{align*}
with the carry operator and the remainder
\begin{align*}
    E_A:=\ket{1^j}\bra{0^j},
    \qquad
    \bar{S}_R:=S_R-\mathds{1}_R ,
\end{align*}
and we write $E_R:=\ket{1^{\,n-1-j}}\bra{0^{\,n-1-j}}$ for the corresponding
operator on the other side of the cut. Because $E_A$ maps a single basis state
to a single basis state, products involving it collapse; the relations needed
below are
\begin{align*}
X_AE_A &= P_A,\qquad
P_AX_AS_A = P_A,\\
E_AX_AS_A &= E_A,\qquad
E_AP_A = E_A, \\
S_AP_A &= E_A.
\end{align*}
Applying them, every $W_{ab}$ reduces from four product terms to two, and the
eight surviving terms use only six distinct left factors. Collecting them, $V$
takes the form \eqref{eq:OSD} with
\begin{equation*}
\begin{aligned}
A_1 &= \mathds{1}_A, &\quad
A_2 &= P_A, &\quad
A_3 &= X_AS_A,\\
A_4 &= S_A, &\quad
A_5 &= E_A, &\quad
A_6 &= S_AX_AS_A,
\end{aligned}
\end{equation*}
and
\begin{equation}
\begin{aligned}
B_1 =\;& b_{00}\,\mathds{1}_R\otimes\ket{0}\bra{0},\\
B_2 =\;& \left(b_{01}X_R\bar{S}_R+d_{01}P_RX_RS_R\right)\otimes\ket{0}\bra{1}\\
&+d_{00}\,P_R\otimes\ket{0}\bra{0},\\
B_3 =\;& b_{01}\,X_R\otimes\ket{0}\bra{1},\\
B_4 =\;& b_{10}\,\mathds{1}_R\otimes\ket{1}\bra{0},\\
B_5 =\;
&\Big[b_{11}\left(X_R\bar{S}_R+\bar{S}_RX_R+\bar{S}_RX_R\bar{S}_R\right)\\
 &+d_{11}\,E_RX_RS_R\Big]\otimes\ket{1}\bra{1}\\
&+\left(b_{10}\bar{S}_R+d_{10}E_R\right)\otimes\ket{1}\bra{0},\\
B_6 =\;& b_{11}\,X_R\otimes\ket{1}\bra{1}.
\end{aligned}
\label{eq:six_right_factors}
\end{equation}
By \eqref{eq:rank_product_bound} this already sharpens
\eqref{eq:loose_estimate} from $36$ to $6$. It remains to show that six is
attained, that is, that the left factors are linearly independent. Their Gram
matrix decouples into two identical blocks, one spanned by
$\{\mathds{1}_A,P_A,X_AS_A\}$ and one by $\{S_A,E_A,S_AX_AS_A\}$, each of the
form
\begin{align*}
G_{\mathrm{block}}
=
\begin{pmatrix}
2^j & 1 & 2\\
1 & 1 & 1\\
2 & 1 & 2^j
\end{pmatrix},
\quad
\det\left(G_{\mathrm{block}}\right)=2^j\big(2^j-2\big),
\end{align*}
so that $\det G=[2^j(2^j-2)]^2$. For $j\ge2$ the determinant is non-zero and
all six factors are independent, whereas $j=1$ makes it vanish and the Gram
matrix drops to rank four. Hence
\begin{align}
    \operatorname{rank}_{S}(V)=
    \begin{cases}
        4, & j=1\ \text{or}\ j=n-1,\\
        6, & 2\leq j\leq n-2 ,
    \end{cases}
    \label{eq:V_exact_rank}
\end{align}
so the factor six is an exact property of the right-hand block rather than an
artefact of the product bound.

\subsection{Boundary-fixed circuit and the expressivity bound}
\label{app:boundary_fixed}

For the SFCT and SFST the physical index of the final qubit is fixed, so only
half of the terms in \eqref{eq:six_right_factors} contribute. For a fixed
output index $a_0$ the surviving left factors are
$\{\mathds{1}_A,P_A,X_AS_A\}$ for $a_0=0$ and
$\{S_A,E_A,S_AX_AS_A\}$ for $a_0=1$ -- exactly the two Gram blocks above, each
of full rank for $j\ge2$. Hence
$\operatorname{rank}_{S}(V_{\mathrm{fixed}})=3$ for $2\leq j\leq n-2$,
reducing to two at the outermost cut.

Combining $\operatorname{rank}_{S}(U_1)=2$, the rank $r_S$ of the QFT and
\eqref{eq:V_exact_rank} gives
\begin{align*}
    \hat r\leq2\,r_S\,6=12\,r_S
\end{align*}
for the QCST with open boundary indices, and
\begin{align*}
    \hat r\leq2\,r_S\,3=6\,r_S
\end{align*}
for the boundary-fixed SFCT and SFST. Both bounds are independent of the
discretisation order $n$, which establishes that the SFCT and SFST admit
efficient TT representations with a saturating rank.

One qualification remains. The rank of the assembled circuit can fall below the
product bound, since further linear dependencies arise when the right-hand
block is composed with the QFT. Equation \eqref{eq:V_exact_rank} is therefore a
rigorous circuit-level bound rather than a claim that the full circuit attains
it.

\begin{figure*}[t]
    \centering
    \includegraphics[width=\linewidth]{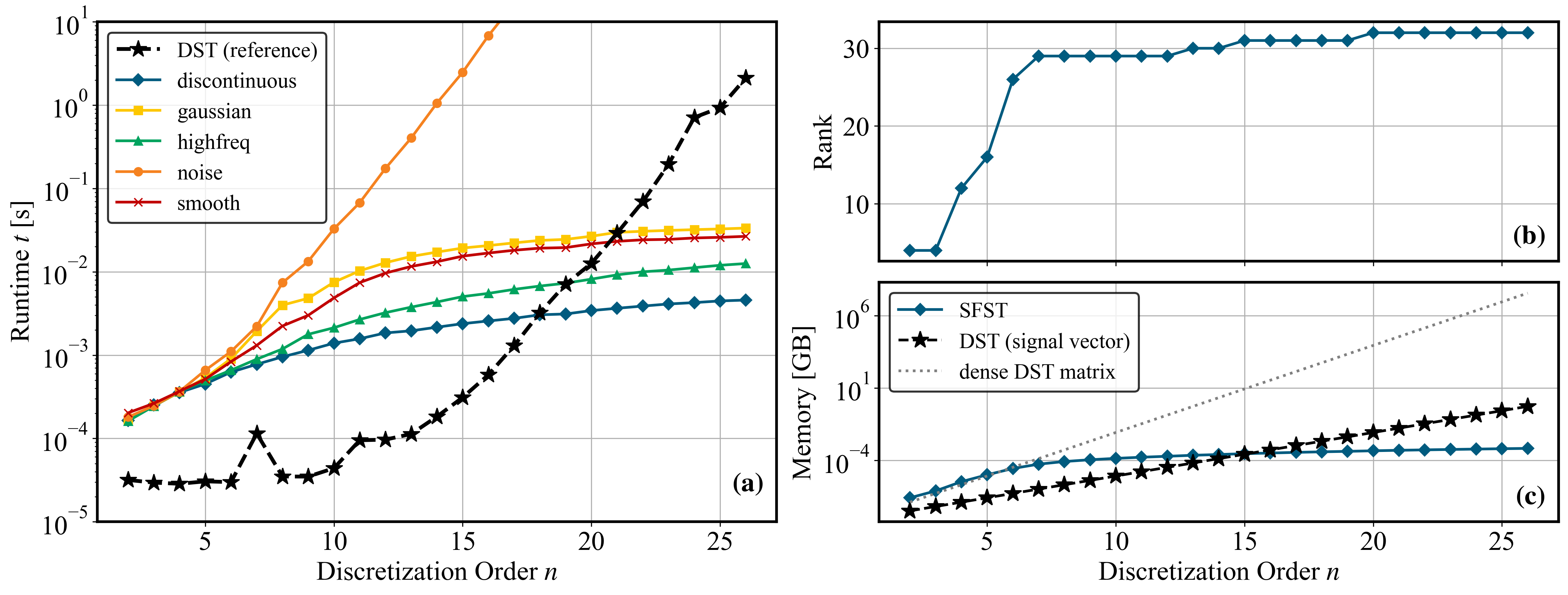}
    \caption{
    (a) Runtime comparison of the classical DST and the SFST for different input classes as a function of the discretization order $n$;
    (b) TTO rank evolution of the SFST;
    (c) Memory of the SFST operator compared with the signal vector required by the classical
    FFT-based DST (dashed) and, for reference, the dense DST matrix (dotted).
    }
    \label{fig:exp_DSTvsSFST}
\end{figure*}

\section{Input Functions}
\label{app:SFCT_Inputs}

We provide the explicit definitions of the input functions used in the numerical experiments of Section~\ref{subsec:numerical_results} and Appendix~\ref{app:exps_dct2}. All functions are defined on $x \in [0,1)$ and discretized on a uniform grid of size $N = 2^n$ via
\begin{align*}
    x_k = \frac{k}{N}, \quad k = 0, \dots, N-1.
\end{align*}

\paragraph{Smooth signal.}
\begin{align*}
    f_{\mathrm{smooth}}(x) = \cos(2\pi x) + \tfrac{1}{2}\cos(4\pi x).
\end{align*}

\paragraph{High-frequency signal.}
\begin{align*}
    f_{\mathrm{high}}(x) = \cos(20\pi x).
\end{align*}

\paragraph{Discontinuous signal.}
\begin{align*}
    f_{\mathrm{disc}}(x) =
    \begin{cases}
        1, & x < 0.5, \\
        0, & x \geq 0.5.
    \end{cases}
\end{align*}

\paragraph{Gaussian signal.}
\begin{align*}
    f_{\mathrm{gauss}}(x) = \exp\big(-100(x - 0.5)^2\big).
\end{align*}

\paragraph{Noise signal.}
We additionally consider a stochastic input given by Gaussian white noise:
\begin{align*}
    f_{\mathrm{noise}}(x_k) \sim \mathcal{N}(0,1),
\end{align*}
with independent samples at each grid point.

These test cases include both structured deterministic signals and a stochastic signal, allowing us to assess the robustness of the SFCT and SFST with respect to different input regularities and noise levels.

\section{Numerical Experiments on DST-II}
\label{app:exps_dct2}
Figure~\ref{fig:exp_DSTvsSFST} presents the results of the experiments
in Section~\ref{subsec:numerical_results} for the DST-II and the SFST using the functions defined in Appendix~\ref{app:SFCT_Inputs}.
The runtime and memory trends are consistent with those observed for the
DCT-II and the SFCT, respectively. Moreover, the TT ranks coincide for
all tested values of $n$, as expected from the equality of the exact
operator Schmidt ranks of the SFCT and SFST across every bipartition.

\section{Preconditioned linear solve in TT format}
\label{app:preconditioned_als}

Poor conditioning of linear systems can slow iterative methods, including TT-based solvers.
In the alternating linear scheme (ALS)~\cite{Holtz2012}, this can
increase both the number of forward--backward sweeps and the effort
required to solve each local linear system.

To remove the constant null mode, we replace the row and column of
$D_h\nabla_h$ associated with one pressure grid point by those of
$-I$, where $I$ is the identity matrix. We denote the resulting matrix
by $A$ and refer to that location as the pinned point.
We represent the pressure $p$ and RHS $b$ as TTs and $A$ and the
preconditioner $P$ as TTOs, following Section~\ref{sec:tensor_trains}.
For the Rayleigh--B\'enard benchmark, we construct a left preconditioner $P$ once before time integration and reuse it in all runs with the same pressure matrix $A$. At each pressure solve, ALS solves the preconditioned system $PAp=Pb$ in TT format~\cite{VANHULST2026110169}.

\paragraph{Offline construction.}
\begin{samepage}
We construct an approximate left inverse $P$ by solving $PA=I$ in TT form.
Let $\operatorname{vec}$ denote column-wise vectorization and
$\otimes$ the Kronecker product. The vectorized system is
\begin{equation}
    (A^{\mathsf T}\otimes I)\operatorname{vec}(P)
    =\operatorname{vec}(I).
    \label{eq:rb_preconditioner_construction}
\end{equation}
\end{samepage}
We use single-site ALS with residual enrichment after each sweep,
inspired by AMEn~\cite{Dolgov2013}, and cap the
operator ranks at $r_{\mathrm P}$.

\begin{figure*}[!t]
    \centering
    \includegraphics[width=.88\textwidth]{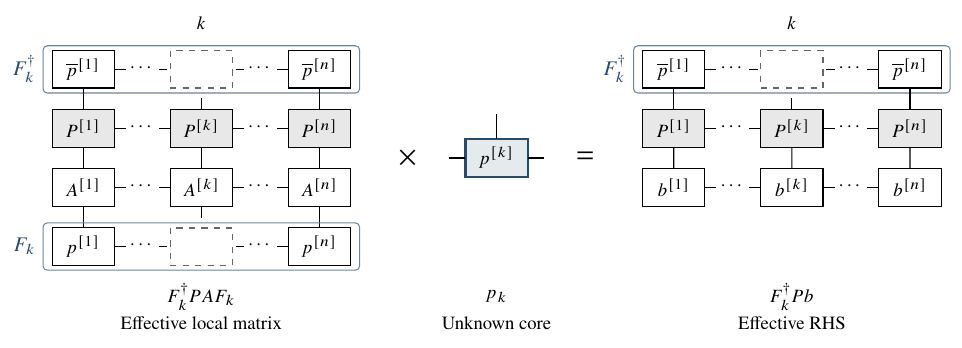}
    \caption{Effective local linear system in
    Eq.~\eqref{eq:rb_preconditioned_local}, with the unknown pressure core
    shown separately. Superscripts $[j]$ label TT or TTO cores, with
    $n=n_x+n_y$, and overbars denote complex conjugation.
    Dashed outlines mark omitted cores, leaving three input and three
    output legs in the effective matrix. All other pressure cores are fixed.}
    \label{fig:rb_als_network}
\end{figure*}

\paragraph{Online pressure solve.}
At site $k$, the unknowns are the entries $p_k=\operatorname{vec}(p^{[k]})$ of the core shown separately in Fig.~\ref{fig:rb_als_network}. The other pressure cores are fixed and define the frame matrix $F_k$~\cite{Dolgov2013}, with $p=F_kp_k$.
Substituting this representation into the preconditioned system and projecting with $F_k^\dagger$
gives the one-site ALS system~\cite{Holtz2012,Schollwock2011}
\begin{equation}
    (F_k^\dagger PAF_k)p_k=F_k^\dagger Pb.
    \label{eq:rb_preconditioned_local}
\end{equation}
Figure~\ref{fig:rb_als_network} shows the contractions with $P$ and $A$
retained as separate TTOs.
We solve this local linear system iteratively with GMRES.
After the sweep, we set the pressure to zero at the pinned point.

\section{Numerical details for Rayleigh--B\'enard convection}
\label{app:rb_comparison}

All spectral runs in Figs.~\ref{fig:rb_observables} and~\ref{fig:rb_matched_accuracy} use $\varepsilon_{\mathrm{spec}}=10^{-10}$.

For the timing comparison in Fig.~\ref{fig:rb_matched_accuracy},
all four configurations advance the same compressed initial state for 12 steps
at $r=100$. The preconditioner rank is $r_{\mathrm P}=50$.
Preconditioned ALS uses one forward--backward sweep with local GMRES
absolute and relative tolerances of $10^{-8}$. VAR uses at most five
forward--backward sweeps, stopping when the relative change in the
objective or iterate falls below $10^{-8}$.

RSI uses an intermediate rank cap $r_{\mathrm{int}}=r+30$.
For the binary TT cores, the local Gaussian sketch matrices have size
$2\times m_{\mathrm{sk}}$~\cite{Meng2026-arxiv},
where $m_{\mathrm{sk}}=\lfloor r_{\mathrm{int}}/2\rfloor+25$ is the sketch
width. RSI starts with a block of eight TT cores and uses a relative
tolerance of $10^{-15}$ for the interpolative decompositions. The product
is subsequently rounded to rank cap $r$.

We perform three sequential repetitions on an AMD EPYC 9654 processor using
Julia 1.12.6 with one Julia thread and one BLAS thread. Method order
varies between repetitions, while the random seeds are held fixed.
The RSI sketch seed is reset at each call.
For each repetition, we average the time per step over steps 3--10.
Pressure-solve time sums both stages and includes convergence tests,
gauge fixing and realification. Full-step time also includes RHS
formation and all state updates. Neither includes setup, diagnostic
reporting or output.
We report the median of the three within-repetition speedup ratios
relative to preconditioned ALS.

As a check on the pressure solves in the timing runs of Fig.~\ref{fig:rb_matched_accuracy}, we evaluated the relative residual of the unpinned pressure equation, $\|b-D_h\nabla_h p\|_2/\|b\|_2$, outside the timed region. Over both projection stages of steps 3--10, the median was $1.04\times10^{-3}$ for preconditioned ALS and ranged from $6.20\times10^{-4}$ to $6.53\times10^{-4}$ across the three spectral configurations.

\section{Real-core TT representation of real-valued states}
\label{app:realification}

A real-valued state can have a TT representation with complex cores.
In both spectral pressure solvers, the inverse Fourier--cosine transform gives a real pressure field whose TT cores may still be complex.
We construct an equivalent representation with real cores by doubling
the internal bond dimensions. For a complex matrix \(M\), define
\begin{equation}
\label{eq:realification_matrix}
\mathscr{R}(M)
=
\begin{pmatrix}
\operatorname{Re} M & -\operatorname{Im} M \\
\operatorname{Im} M &  \operatorname{Re} M
\end{pmatrix}.
\end{equation}
For complex matrices $M_1$ and $M_2$ of compatible dimensions,
$\mathscr{R}(M_1M_2)=\mathscr{R}(M_1)\mathscr{R}(M_2)$,
so this map preserves matrix products.

Write the complex TT as in Eq.~\eqref{eq:tt_format},
\[
T(i_1,\ldots,i_n)
=
A_1(i_1)A_2(i_2)\cdots A_n(i_n),
\]
with open boundary ranks \(r_0=r_n=1\). We apply \(\mathscr{R}\)
to each matrix slice and select the real component at the two boundaries.
With $e_0=(1,0)^{\top}$, the real cores are
\[
\begin{aligned}
B_1(i_1) &= e_0^{\top}\mathscr{R}\!\left(A_1(i_1)\right), \\
B_k(i_k) &= \mathscr{R}\!\left(A_k(i_k)\right),
\qquad k=2,\ldots,n-1, \\
B_n(i_n) &= \mathscr{R}\!\left(A_n(i_n)\right)e_0 .
\end{aligned}
\]
The contraction of the realified TT gives
\begin{equation}
\begin{aligned}
T_{\mathbb{R}}(i_1,\ldots,i_n)
&=
B_1(i_1)B_2(i_2)\cdots B_n(i_n) \\
&=
e_0^{\top}
\mathscr{R}\!\left(
A_1(i_1)A_2(i_2)\cdots A_n(i_n)
\right)
e_0 \\
&=
\operatorname{Re}T(i_1,\ldots,i_n).
\end{aligned}
\end{equation}
For a real-valued tensor, the construction preserves the represented
values up to roundoff.
Before recompression, the internal bond dimensions are
$\tilde r_k=2r_k$ for $k=1,\ldots,n-1$, with
$\tilde r_0=\tilde r_n=1$. In the pressure solver, the realified
TT is rounded to rank cap $r$.

\section{Validation of Homogenization}
\label{app:validation}

\begin{figure}
    \centering
    \includegraphics[width=\linewidth]{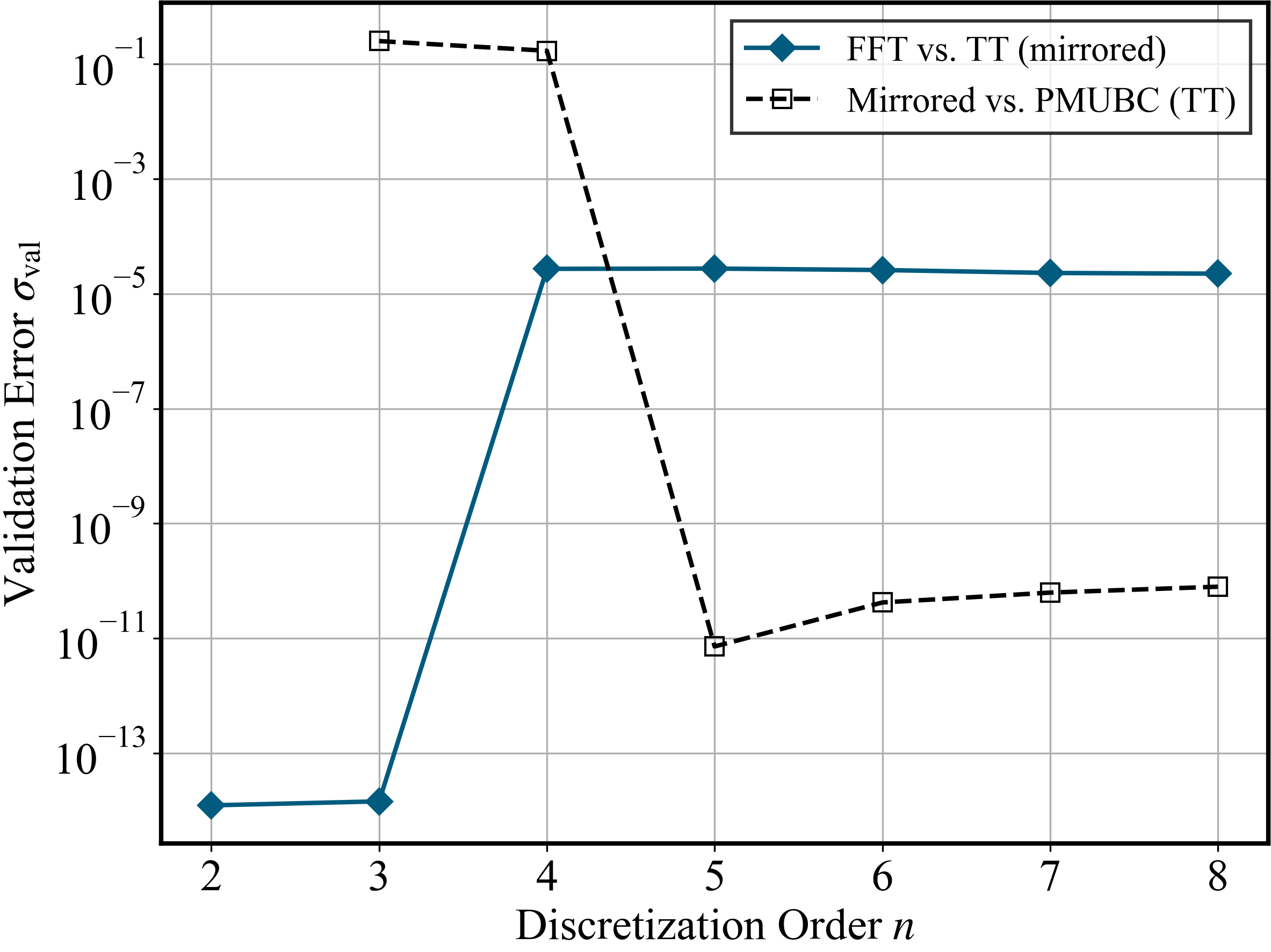}
    \caption{
    Validation of the Spectral-TT Homogenization method by comparing TT-based and FFT-based homogenization using the mirrored-domain approach~(solid blue), as well as the mirrored-domain approach and direct PMUBC using TT-based homogenization~(dashed black).
    The relative error in the homogenized stiffness tensor is shown as a function of the discretization order $n$, with $2^n$ grid points per dimension.
    }
    \label{fig:val_homogenization}
\end{figure}

The figure of merit for the homogenization procedure is the macroscopic stiffness tensor obtained from the Fourier-based algorithms. We assess the accuracy of the TT-based homogenization method through two complementary comparisons, which are shown in Figure~\ref{fig:val_homogenization}.

First, we compare the stiffness tensors obtained using TT-based and conventional FFT-based homogenization, both employing the mirrored-domain representation shown in Figure~\ref{fig:mirrored_domain}. The corresponding relative validation error is defined as
\begin{align*}
    \sigma_\mathrm{val}^\mathrm{(F)} =
    \frac{\lVert\mathbf{C}_\mathrm{TT} - \mathbf{C}_\mathrm{FFT}\rVert_\mathrm{F}}
    {\lVert\mathbf{C}_\mathrm{FFT}\rVert_\mathrm{F}}.
\end{align*}
This comparison serves to validate the TT implementation against the established FFT-based approach over the range of discretization orders for which the conventional FFT method remains computationally feasible.

For the lowest two discretization orders, $n=2$ and $n=3$, the relative error is below $10^{-13}$, indicating almost identical results. For higher orders, the error increases slightly above $10^{-5}$ .
This behavior is consistent with the accuracy controlled by the prescribed truncation tolerance $\delta_\mathrm{tol}=10^{-4}$ in the homogenization algorithms.
Thus, the observed differences are within the expected accuracy of the TT-based procedure and could be further reduced by reducing $\delta_\mathrm{tol}$.

Second, we investigate the influence of the domain representation by comparing the homogenized stiffness tensors obtained using the mirrored-domain approach and direct PMUBC, with both calculations employing the TT-based homogenization algorithm. The corresponding relative difference is defined as
\begin{align*}
    \sigma_\mathrm{val}^\mathrm{(P)} =
    \frac{\lVert\mathbf{C}_\mathrm{mirrored} - \mathbf{C}_\mathrm{PMUBC}\rVert_\mathrm{F}}
    {\lVert\mathbf{C}_\mathrm{PMUBC}\rVert_\mathrm{F}}.
\end{align*}

For $n=3$ and $n=4$, the relative difference is initially on the order of $10^{-1}$. With increasing discretization order, however, it decreases rapidly to approximately $10^{-11}$ and subsequently remains at a level of around $10^{-10}$. This indicates that, for sufficiently fine discretizations, both domain representations yield virtually identical homogenized stiffness tensors.

\section{Symmetries and Spectral Methods}
\label{app:PMUBC}
The present work employs
the Hex8R discretization~\cite{Schneider2017}, i.e., trilinear hexahedral elements with reduced integration on a regular voxel grid.
The displacements are located at the voxel corners, whereas
all strain and stress components are evaluated at the single
integration point in the voxel center, see Table~\ref{tab:tab1}. The FFT-based scheme transforms only strain-type quantities. For
these, every reflection symmetry of the mirrored domain is
therefore half-sample symmetric with $f_j=\pm f_{N-1-j}$ and the
required transforms are exclusively the DCT-II for even and
the DST-II for odd symmetry.
Table~\ref{tab:tab2} lists the symmetry of each strain component
with respect to the three coordinate directions for all six unit
load cases.
Each letter selects the transform along the
corresponding direction: E denotes a DCT-II and O a DSTII. For example, under tensile load in $x$-direction, the shear
strain $\varepsilon_{12}$ with pattern OOE is transformed by a DST-II in $x$
and $y$ and by a DCT-II in $z$.

\begin{table}[h]
\centering
\caption{Location of displacement and strain components for Hex8R discretization~\cite{Schneider2017}.}
\label{tab:tab1}
\begin{tabular}{c c}
\hline
Variable & Location \\
\hline
$u_1, u_2, u_3$ & $(i h_1, j h_2, k h_3)$ \\
$\varepsilon_{11}, \varepsilon_{22}, \varepsilon_{33}$ & $\left(\left(i+\tfrac{1}{2}\right)h_1,\left(j+\tfrac{1}{2}\right)h_2,\left(k+\tfrac{1}{2}\right)h_3\right)$ \\
$\varepsilon_{12}, \varepsilon_{13}, \varepsilon_{23}$ & $\left(\left(i+\tfrac{1}{2}\right)h_1,\left(j+\tfrac{1}{2}\right)h_2,\left(k+\tfrac{1}{2}\right)h_3\right)$ \\
\hline
\end{tabular}
\end{table}

\begin{table}[h]
\centering
\caption{Symmetry of strain components for different loading cases.
E/O denote even/odd symmetry with respect to $(x,y,z)$~\cite{Grimm2021}.}
\label{tab:tab2}
\begin{tabular}{c c c c c c c}
\hline
Load & $\varepsilon_{11}$ & $\varepsilon_{22}$ & $\varepsilon_{33}$ & $\varepsilon_{23}$ & $\varepsilon_{13}$ & $\varepsilon_{12}$ \\
\hline
$11$ & EEE & EEE & EEE & EOO & OEO & OOE \\
$22$ & EEE & EEE & EEE & EOO & OEO & OOE \\
$33$ & EEE & EEE & EEE & EOO & OEO & OOE \\
$23$ & EOO & EOO & EOO & EEE & OOE & OEO \\
$13$ & OEO & OEO & OEO & OOE & EEE & EOO \\
$12$ & OOE & OOE & OOE & OEO & EOO & EEE \\
\hline
\end{tabular}
\end{table}

\printcredits

\bibliographystyle{model1-num-names}

\bibliography{02_main}

\end{document}